\documentclass[10pt]{article}
\usepackage{amsfonts,amssymb,amsmath,comment}
\usepackage{graphicx,graphics,color}
\usepackage[T1]{fontenc}
\usepackage[utf8]{inputenc}
\usepackage{graphics,epsfig,color}
\usepackage{setspace,enumitem,verbatim}
\usepackage{caption}
\usepackage{subcaption}
\usepackage{etoolbox}
\usepackage{enumerate}

 \makeatletter
 \@addtoreset{equation}{section}
 \makeatother

 \newcounter{enunciato}[section]

 \renewcommand{\theequation}{\thesection.\arabic{equation}}

 \newtheorem{ittheorem}{Theorem}
 \newtheorem{itlemma}{Lemma}
 \newtheorem{itproposition}{Proposition}
 \newtheorem{itdefinition}{Definition}
 \newtheorem{itcorollary}{Corollary}
 \newtheorem{itconjecture}{Conjecture}
 \newtheorem{itremark}{Remark}
 \newtheorem{itclaim}{Claim}
\newtheorem{assumption}{Assumption}

 \newenvironment{theorem}{\addtocounter{enunciato}{1}
 \begin{ittheorem}}{\end{ittheorem}}

 \newenvironment{lemma}{\addtocounter{enunciato}{1}
 \begin{itlemma}}{\end{itlemma}}

 \newenvironment{proposition}{\addtocounter{enunciato}{1}
 \begin{itproposition}}{\end{itproposition}}

 \newenvironment{definition}{\addtocounter{enunciato}{1}
 \begin{itdefinition}}{\end{itdefinition}}

 \newenvironment{remark}{\addtocounter{enunciato}{1}
 \begin{itremark}}{\end{itremark}}

\newcommand{\halmos}{\rule{1ex}{1.4ex}}

\newenvironment{proof}{\noindent {\em Proof}.\,\,}
{\hspace*{\fill}$\halmos$\medskip}

\newcommand{\N}{\mathbb{N}}

\newcommand{\cG}{\mathcal{G}}

\newcommand{\vC}{\vec{C}}
\newcommand{\vt}{\vec{\theta}}
\newcommand{\Pmic}{\mathrm{P}_{\mathrm{mic}}}
\newcommand{\Pcan}{\mathrm{P}_{\mathrm{can}}}
\newcommand{\dd}{\mathrm{d}}
\newcommand{\e}{\epsilon}
\newcommand{\eee}{\mathrm{e}}

\newcommand\numberthis{\addtocounter{equation}{1}\tag{\theequation}}

\allowdisplaybreaks

\begin{document}

\title{Breaking of ensemble equivalence for\\ 
perturbed Erd\H{o}s-R\'enyi random graphs}

\author{\renewcommand{\thefootnote}{\arabic{footnote}}
F. den Hollander\footnotemark[1]\,,
\renewcommand{\thefootnote}{\arabic{footnote}}
M. Mandjes\footnotemark[2]\,,
\renewcommand{\thefootnote}{\arabic{footnote}}
A. Roccaverde\footnotemark[3]\,,
\renewcommand{\thefootnote}{\arabic{footnote}}
N.J. Starreveld\footnotemark[4]
}

\footnotetext[1]{Mathematical Institute, Leiden University, P.O. Box 9512, 2300 RA Leiden,
The Netherlands\\ 
{\tt \small denholla@math.leidenuniv.nl}}

\footnotetext[2]{Korteweg de-Vries Institute, University of Amsterdam, P.O. Box 94248, 
1090 GE Amsterdam, The Netherlands\\ 
{\tt \small m.r.h.mandjes@uva.nl}}

\footnotetext[3]{Mathematical Institute, Leiden University, P.O. Box 9512, 2300 RA Leiden, 
The Netherlands\\ 
{\tt \small roccaverdeandrea@gmail.com }}

\footnotetext[4]{Korteweg de-Vries Institute, University of Amsterdam, P.O. Box 94248, 
1090 GE Amsterdam, The Netherlands\\ 
{\tt \small n.j.starreveld@uva.nl}}

\date{\today}

\maketitle


\begin{abstract}

In \cite{dHMRS18} we analysed a simple undirected random graph subject to constraints on the 
densities of edges and triangles, considering the dense regime in which the number of edges 
per vertex is proportional to the number of vertices. We computed the specific relative entropy of 
the \emph{microcanonical ensemble} with respect to the \emph{canonical ensemble}, i.e., the 
relative entropy per edge in the limit as the number of vertices tends to infinity. We showed that 
as soon as the constraints are \emph{frustrated}, i.e., do not lie on the Erd\H{o}s-R\'enyi line (where
the density of triangles is the third power of the density of edges), there is \emph{breaking of ensemble 
equivalence}, meaning that the specific relative entropy is strictly positive. In the present paper 
we analyse what happens near this line. It turns out that the way in which the specific relative 
entropy tends to zero critically depends on whether the line is approached is from above or from below. We 
identify what the constrained random graph looks like in the microcanonical ensemble in the limit 
as the number of vertices tends to infinity. 
 
\vspace{0.5cm}\noindent
\emph{MSC 2010:} 05C80, 60K35, 82B20.\\ 
\emph{Key words:} Erd\H{o}s-R\'enyi random graph, Gibbs ensembles, breaking of ensemble equivalence,
relative entropy, graphons, large deviation principle, variational representation.\\
\emph{Acknowledgement:} The research in this paper was supported through NWO Gravitation Grant 
NETWORKS 024.002.003. The authors are grateful to V.\ Patel and H.\ Touchette for helpful discussions. 
FdH, AR and NJS are grateful for hospitality at the International Centre for Theoretical Sciences in 
Bangalore, India, as participants of the program on \emph{Large Deviation Theory in Statistical Physics: 
Recent Advances and Future Challenges} running in the Fall of 2017.

\end{abstract}


\newpage




\section{Introduction}
\label{S1.1}


\subsection{Background}

In this paper we analyse random graphs that are subject to \emph{constraints}. Statistical 
physics prescribes what probability distribution on the set of graphs we should choose 
when we want to model a given type of constraint \cite{G02}. Two important choices are: 
\begin{itemize}
\item[(1)] 
The \emph{microcanonical ensemble}, where the constraints are \emph{hard} (i.e., are 
satisfied by each individual graph).
\item[(2)] 
The \emph{canonical ensemble}, where the constraints are \emph{soft} (i.e., hold as 
ensemble averages, while individual graphs may violate the constraints).
\end{itemize}
For random graphs that are large but finite, the two ensembles are obviously different and, 
in fact, represent different empirical situations. Each ensemble represents the unique 
probability distribution with \emph{maximal entropy} respecting the constraints. In the 
limit as the size of the graph diverges, the two ensembles are traditionally \emph{assumed} 
to become equivalent as a result of the vanishing fluctuations in the soft constraints, i.e., 
the soft constraints are assumed to behave asymptotically like hard constraints. This 
assumption of \emph{ensemble equivalence} is one of the cornerstones of statistical physics, 
but it does \emph{not} hold in general. We refer to \cite{T14} for more background on this 
phenomenon.

In a series of papers we investigated the possible breaking of ensemble equivalence for 
various choices of the constraints, including the degree sequence and the total number of 
edges, wedges, triangles, etc. Both the \emph{sparse regime} (where the number of edges 
per vertex remains bounded) and the \emph{dense regime} (where the number of edges per 
vertex is of the order of the number of vertices) were considered. The effect of \emph{community 
structure} on ensemble equivalence was investigated as well. Relevant references are \cite{GHR17,GHR18,dHMRS18,GS18,SdMdHG15}. 

In \cite{dHMRS18}, for the dense regime, we considered constraints on the densities of 
\emph{finitely many arbitrary subgraphs}. Our main result was a \emph{variational formula} for 
$s_\infty = \lim_{n\to\infty} n^{-2} s_n$, where $n$ is the number of vertices and $s_n$ is the 
relative entropy of the microcanonical ensemble with respect to the canonical ensemble. 
Our analysis relied on the \emph{large deviation principle for graphons} derived in \cite{Ch15,CV11}. 
For the case where the constraints were on the densities of edges and triangles we found that 
$s_\infty>0$ when the constraints are \emph{frustrated}.

In the sequel we will say that we are on the {\em Erd\H{o}s-R\'enyi line} (abbreviated to ER 
line) when the density of triangles is equal to the third power of the density of edges, which is
typical for the Erd\H{o}s-R\'enyi random graph. We analyse the behaviour of $s_\infty$ when 
the constraints are close to but different from the ER line. Moreover, we identify what the 
constrained random graph looks like asymptotically in the microcanonical ensemble. It turns 
out that the behaviour drastically changes when the density of triangles is slightly larger, 
respectively, slightly smaller than that of the ER line. The microcanonical ensemble is harder 
to analyse than the canonical ensemble. Yet, we do not know what the constrained random 
graph looks like asymptotically in the canonical ensemble below the ER line.  


\subsection{Literature}

While breaking of ensemble equivalence is a relatively new concept in the theory of random graphs, 
there are many studies on the asymptotic structure of random graphs. In the pioneering work \cite{CV11}, 
followed by \cite{LZ15}, a large deviation principle for dense Erd\H{o}s-R\'enyi random graphs was
established and the asymptotic structure of constrained Erd\H{o}s-R\'enyi random graphs was 
described as the solution of a variational problem. In the past few years, significant progress was 
made regarding sparse random graphs as well \cite{CDS11,DL18,LZ16,Z17}. Two further random 
graph models that were studied extensively are the exponential random graph model and the 
constrained exponential random graph model. Exponential random graphs, which are related to 
the canonical ensemble considered in the present paper, were analysed in \cite{BBS11,CD13}: \cite{BBS11} 
studies mixing times of Glauber spin-flip dynamics, while \cite{CD13} uses large deviation theory to 
derive asymptotic expressions for the partition function. In subsequent works \cite{RS13,RY13,Y13} 
the behaviour of exponential random graphs was analysed in further detail, while in \cite{YZ17} the focus 
was on sparse exponential random graphs. In \cite{EG18} exponential random graphs were studied 
in both the dense and the sparse regime, and the main conclusion was that they behave essentially 
like mixtures of random graphs with independent edges. In \cite{KY17,Y15} constrained exponential 
random graphs were analysed, while in \cite{AZ11,AZ18} the additional feature of directed edges was 
added. In \cite{DS19} large deviations were used to study random graphs constrained on the 
degree sequence in the dense regime. 

In \cite{KRRS17,KRRS172,MY18,RS15}, the asymptotic structure of graphs drawn from the microcanonical 
ensemble was investigated for various choices of constraints on the densities of edges and triangles. 
The focus of \cite{RS15} was on the behaviour of random graphs for values of the edge and triangle 
densities close to the ER line, and a rough scaling of the graph was found via a bound on the entropy 
function. In the present paper we extend these results by determining the precise scaling. A similar 
question was addressed in \cite{MY18} for a constraint on the edge and triangle density close to the 
lower boundary curve of the admissibility region (see Fig.~\ref{fig-scallopy} below). In \cite{KRRS172}, 
through extensive numerics, the regions were determined where phase transitions in the structure 
of the constrained random graph occur as the densities of edges and triangles is varied. Our results 
rely on this numerics near the ER line, and make it mathematically precise.  


\subsection{Outline}

The remainder of our paper is organised as follows. In Section~\ref{S1.2} we define the two ensembles, 
give the definition of equivalence of ensembles in the dense regime, and recall some basic facts about 
graphons. An important role is played by the \emph{variational representation} of $s_\infty$, derived in 
\cite{dHMRS18} when the constraints are on the total numbers of subgraphs drawn from a finite collection 
of subgraphs. We also recall the analysis of $s_\infty$ in \cite{dHMRS18} for the special case where the 
subgraphs are the edges and the triangles. In Section~\ref{S1.4} we state our main theorems and 
propositions on the behaviour around the ER line. Proofs are given in Sections~\ref{S2 Chapter 5} and 
\ref{Optimal per}.



\section{Definitions and preliminaries}
\label{S1.2}

In this section, which is largely lifted from \cite{dHMRS18}, we present the definitions of the main concepts 
to be used in the sequel, together with some key results from prior work. Section \ref{S1.2.1} presents the 
formal definition of the two ensembles we are interested in and gives our definition of ensemble equivalence  
in the dense regime. Section \ref{S1.2.2} recalls some basic facts about  {graphons}, while Section \ref{S1.2.3} 
recalls some basic properties of the canonical ensemble. Section \ref{S1.3} recalls the variational characterisation 
of ensemble equivalence when the constraint is on a finite number of subgraph densities, proven in \cite{dHMRS18}. Section \ref{Sedtr} recalls the main results in \cite{dHMRS18} for the case where the constraints are on the densities 
of edges and triangles.   


\subsection{Microcanonical ensemble, canonical ensemble, relative entropy}
\label{S1.2.1}

For $n \in \N$, let $\cG_n$ denote the set of all $2^{{n\choose 2}}$ simple undirected graphs with vertex 
set $\{1,\dots, n\}$. Any graph $G\in\cG_n$ can be represented by a symmetric $n \times n$ matrix 
with elements 
\begin{equation}
h^G(i,j) :=
\begin{cases}
1\qquad \mbox{if there is an edge between vertex } i \mbox{ and vertex } j,\\ 
0 \qquad \mbox{otherwise.}
\end{cases}
\end{equation}
Let $\vC$ denote a vector-valued function on $\cG_n$. We choose a specific vector $\vC^*$,
which we assume to be \emph{graphical}, i.e., realisable by at least one graph in $\cG_n$.
Given $\vC^*$, the \emph{microcanonical ensemble} is the probability distribution $\Pmic$ 
on $\cG_n$ with \emph{hard constraint} $\vC^*$ defined as
\begin{equation}
\Pmic(G) :=
\left\{
\begin{array}{ll} 
1/\Omega_{\vC^*}, \quad & \text{if } \vC(G) = \vC^*, \\ 
0, & \text{otherwise},
\end{array}
\right. \qquad G\in \cG_n,
\label{eq:PM}
\end{equation}
where 
\begin{equation}
\Omega_{\vC^*} := | \{G \in \cG_n\colon\, \vC(G) = \vC^* \} |
\end{equation}
is the number of graphs that realise $\vC^*$. The \emph{canonical ensemble} $\Pcan$ is 
the unique probability distribution on $\cG_n$ that maximises the \emph{entropy} 
\begin{equation}
S_n({\rm P}) := - \sum_{G \in \cG_n}{\rm P}(G) \log {\rm P}(G)
\end{equation}
subject to the \emph{soft constraint} $\langle \vC \rangle  = \vC^*$, where
\begin{equation}
\label{softconstr}
\langle \vC \rangle  := \sum_{G \in \cG_n} \vC(G)\,{\rm P}(G).
\end{equation}
This gives the formula \cite{J57}
\begin{equation}
\Pcan(G) := \frac{1}{Z(\vt^*)}\,\eee^{H(\vec{\theta}^*,\vec{C}(G))},
\qquad G \in \cG_n,
\label{eq:PC}
\end{equation}
with 
\begin{equation}
H(\vec{\theta}^*,\vec{C}(G)) := \vt^* \cdot \vC(G), \qquad
Z(\vt^*\,) := \sum_{G \in \cG_n} \eee^{\vt^* \cdot\hspace{2pt} \vC(G)},
\label{eq:Ham}
\end{equation}
denoting the \emph{Hamiltonian} and the \emph{partition function}, respectively. In 
\eqref{eq:PC}--\eqref{eq:Ham} the parameter $\vt^*$, which is a real-valued vector 
whose dimension is equal to the number of constraints, must be set to the unique value 
that realises $\langle \vC \rangle  = \vC^*$. As a Lagrange multiplier, $\vec{\theta}^*$ 
always exists, but uniqueness is non-trivial. In the sequel we will only consider examples
where the gradients of the constraints are \emph{linearly independent} vectors. Consequently, 
the Hessian matrix of the covariances in the canonical ensemble is a positive-definite matrix, 
which implies uniqueness. 

The \emph{relative entropy} of $\Pmic$ with respect to  $\Pcan$ is defined as
\begin{equation}
S_n(\Pmic \mid \Pcan) 
:= \sum_{G \in \cG_n} \Pmic(G) \log \frac{\Pmic(G)}{\Pcan(G)}.
\label{eq:KL1}
\end{equation}
For any $G_1,G_2\in\cG_n$, $\Pcan(G_1)=\Pcan(G_2)$ whenever $\vC(G_1)=\vC(G_2)$, 
i.e., the canonical probability is the same for all graphs with the same value of the constraint. 
We may therefore rewrite \eqref{eq:KL1} as
\begin{equation}
S_n(\Pmic \mid \Pcan) = \log \frac{\Pmic(G^*)}{\Pcan(G^*)},
\label{eq:KL2}
\end{equation}
where $G^*$ is \emph{any} graph in $\cG_n$ such that $\vC(G^*) =\vC^*$ (recall that we assumed 
that $\vC^*$ is realisable by at least one graph in $\cG_n$). The removal of the sum over $\cG_n$
constitutes a major simplification. 

All the quantities above depend on $n$. In order not to burden the notation, we exhibit this $n$-dependence 
only in the symbols $\cG_n$ and $S_n(\Pmic \mid \Pcan)$. When we pass to the limit $n\to\infty$, we need 
to specify how $\vC(G)$, $\vC^*$ are chosen to depend on $n$. We refer the reader to \cite[Appendix A]{dHMRS18}, 
where this issue was discussed in detail. 

\begin{definition}
\label{def: sinfty}
{\rm In the dense regime, if 
\begin{equation}
\label{eq: def relative entropy limit}
s_{\infty} : = \lim\frac{1}{n^2}S_n(\Pmic \mid \Pcan)=0,
\end{equation}
then $\Pmic$ and $\Pcan$ are said to be {\em equivalent}.}
\end{definition}

\begin{remark}
{\rm In \cite{SdMdHG15}, which was concerned with the \emph{sparse regime}, the relative entropy 
was divided by $n$ (the number of vertices). In the \emph{dense regime}, however, it is appropriate 
to divide by $n^2$ (the order of the number of edges). }
\end{remark} 


\subsection{Graphons}
\label{S1.2.2}

\begin{figure}[htbp]
\centering
\hspace{1.5cm}\includegraphics[width=0.7\linewidth, height=4.6cm]{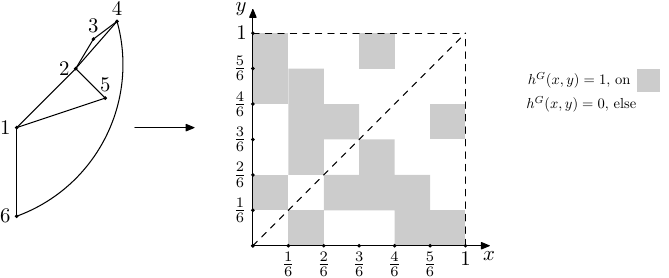}
\caption{\small An example of a graph $G$ and its graphon representation $h^G$.} 
\label{fig-graphon}
\end{figure}

There is a natural way to embed a simple graph on $n$ vertices in a space of functions called 
\emph{graphons}. Let $W$ be the space of functions $h\colon\,[0,1]^2 \to [0,1]$ such that 
$h(x,y) = h(y,x)$ for all $(x,y) \in [0,1]^2$. A finite simple graph $G$ on $n$ vertices can be
represented as a graphon $h^{G} \in W$ in a natural way as (see Figure~\ref{fig-graphon})
\begin{equation}
\label{graphondef}
h^{G}(x,y) := \left\{ \begin{array}{ll} 
1 &\mbox{if there is an edge between vertex } \lceil{nx}\rceil \mbox{ and vertex } \lceil{ny}\rceil,\\
0 &\mbox{otherwise},
\end{array} 
\right.
\end{equation}
which is referred to as the empirical graphon associated with $G$.

\medskip\noindent
The space of graphons $W$ is endowed with the \emph{cut distance}
\begin{equation}
d_{\square} (h_1,h_2) := \sup_{S,T\subset [0,1]} 
\left|\int_{S\times T} \dd x\,\dd y\,[h_1(x,y) - h_2(x,y)]\right|,
\qquad h_1,h_2 \in W.
\end{equation}
On $W$ there is a natural equivalence relation $\equiv$. Let $\Sigma$ be the space of 
measure-preserving bijections $\sigma\colon\, [0,1] \to [0,1]$. Then $h_1(x,y)\equiv h_2(x,y)$ 
if $h_1(x,y) = h_2(\sigma x, \sigma y)$ for some $\sigma\in\Sigma$. This equivalence relation 
yields the quotient space $(\tilde{W},\delta_{\square})$, where $\delta_{\square}$ is the 
metric defined by 
\begin{equation}
\label{deltam}
\delta_{\square}(\tilde{h}_1,\tilde{h}_2) 
:= \inf _{\sigma_1,\sigma_2 \in \Sigma} d_{\square}(h_1^{\sigma_1}, h_2^{\sigma_2}),
\qquad \tilde{h}_1,\tilde{h}_2 \in \tilde{W}.
\end{equation}
For a more detailed description of the structure of the space $(\tilde{W},\delta_{\square})$ we refer 
to \cite{BCLSV08,BCLSV12,DGKR15}. In the sequel we will deal with constraints on the edge and 
triangle density. In the space $W$ the edge and triangle densities of a graphon $h$ are defined by 
\begin{equation}
\label{eq: densityalt}
T_1(h) := \int_{[0,1]^2} \dd x_1 \dd x_2\, h(x_1,x_2), \quad 
T_2(h):= \int_{[0,1]^3} \dd x_1 \dd x_2\dd x_3\, h(x_1,x_2)h(x_2,x_3)h(x_3,x_1).
\end{equation} 
For an element $\tilde{h}$ of the quotient space $\tilde{W}$ we define the edge and triangle density by 
\begin{equation}
T_1(\tilde{h}) = T_1(h), \qquad T_2(\tilde{h}) = T_2(h),
\end{equation}
where $h$ is any representative element of the equivalence class $\tilde{h}$. 
 

\subsection{Subgraph counts}
\label{S1.2.3}

We can label the simple graphs in any order: $\{F_k\}_{k\in\N}$. Let $C_k(G)$ denote the number of 
subgraphs $F_k$ in $G$. In the dense regime, $C_k(G)$ grows like $n^{V_k}$ as $n\to\infty$, where 
$V_k:=|V(F_k)|$ is the number of vertices in $F_k$. For $m \in \N$, consider the following \emph{scaled 
vector-valued function} on $\cG_n$:
\begin{equation}
\vec{C}(G) := \left(\frac{p(F_k)\,C_{k}(G)}{n^{V_k-2}}\right)_{k=1}^m
= n^2\left(\frac{p(F_k)\,C_{k}(G)}{n^{V_k}}\right)_{k=1}^m.
\end{equation}
The term $p(F_k)$ counts the edge-preserving permutations of the vertices of $F_k$ (e.g.\ $2$ for an 
edge, $2$ for a wedge, $6$ for a triangle). The term $C_k(G)/n^{V_k}$ represents the density of $F_k$ 
in $G$. The additional $n^2$ guarantees that the full vector scales like $n^2$, in line with the scaling 
of the large deviation principle for graphons in the Erd\H{o}s-R\'enyi random graph derived in \cite{CV11}. 
For a simple graph $F_k$, let $\text{hom}(F_k,G)$ be the number of homomorphisms from $F_k$ to 
$G$, and define the \emph{homomorphism density} as
\begin{equation}
t(F_k,G) := \frac{\text{hom}(F_k,G)}{n^{V_k}}
= \frac{p(F_k)\,C_k(G)}{n^{V_k}},
\end{equation}
which does not distinguish between permutations of the vertices. In terms of this quantity, 
the Hamiltonian becomes
\begin{equation}
\label{eq:HF}
H(\vec{\theta}, \vec{T}(G))=n^2 \sum_{k=1}^m \theta_k \,t(F_k,G) 
= n^2 (\vec{\theta}\cdot\vec{T}(G)), \qquad G \in \cG_n,
\end{equation}
where 
\begin{equation}
\label{operator}
\vec{T}(G) := \left(t(F_k, G)\right)_{k=1}^m.
\end{equation}
The canonical ensemble with parameter $\vec{\theta}$ thus takes the form
\begin{equation}
\label{eq:CPD}
\Pcan(G \mid \vec{\theta}\,) := \eee^{n^2\big[\vec{\theta}\cdot\vec{T}(G)-\psi_n(\vec{\theta}\,)\big]},
\qquad G \in \cG_n,
\end{equation}
where $\psi_n$ replaces the \emph{partition function} $Z(\vt)$: 
\begin{equation}
\label{eq:PF}
\psi_n(\vec{\theta}) := \frac{1}{n^2}\log\sum_{G\in\cG_n} 
\eee^{n^2 (\vec{\theta}\hspace{2pt}\cdot\hspace{2pt}\vec{T}(G))} = \frac{1}{n^2}\log Z(\vt).
\end{equation}
In the sequel we take $\vec{\theta}$ equal to a specific value $\vec{\theta}^*$ so as to meet 
the soft constraint, i.e.,
\begin{equation}
\label{softconstraint2}
\langle \vec{T} \rangle   = \sum_{G\in\cG_n}\vec{T}(G)\,\Pcan(G) = \vec{T}^*.
\end{equation}
With this choice, the canonical probability becomes
\begin{equation}
\label{canonicprob}
\Pcan(G) = \Pcan(G\mid \vec{\theta^*}).
\end{equation}

Both the constraint $\vec{T}^*$ and the \emph{Lagrange multiplier} $\vec{\theta}^*$ in general depend 
on $n$, i.e., $\vec{T}^*=\vec{T}^*_n$ and $\vec{\theta}^* = \vec{\theta}^*_n$. We consider constraints that 
converge when we pass to the limit $n\to\infty$, i.e., 
\begin{equation}
\label{eq: Assumption T}
\lim_{n\to\infty} \vec{T}^*_n = \vec{T}^*_\infty.
\end{equation}
Consequently, we expect that  
\begin{equation}
\label{eq:Assumption}
\lim_{n\to\infty} \vec{\theta}^*_n = \vec{\theta}^*_\infty.
\end{equation}
In the remainder of this paper we \emph{assume} that \eqref{eq:Assumption} holds.\ If convergence fails, 
then we may still consider subsequential convergence. The subtleties concerning \eqref{eq:Assumption} 
were discussed in detail in \cite[Appendix A]{dHMRS18}.


\subsection{Variational characterisation of ensemble equivalence}
\label{S1.3}

For a graphon $h \in W$ and a simple graph $F$ with vertex set $V(F)$ and edge set $E(F)$, define 
\begin{equation}
t(F,h) := \int_{[0,1]^{|V(F)|}} \dd x_1 \dots \dd x_{|V(F)|.} \prod_{\{i,j\} \in E(F)} h(x_i,x_j).    
\end{equation}
Then $t(F,G) = t(F,h^G)$, and so the  expression in \eqref{eq:HF} can be written as
\begin{equation}
\label{eq:HFG}
H(\vec{\theta}, \vec{T}(G))=n^2 \sum_{k=1}^m \theta_k \,t(F_k,h^G).
\end{equation}
The constraint $\vec{T}^*_\infty$ defines a subspace of the quotient space $\tilde{W}$, 
\begin{equation}
\tilde{W}^* := \{\tilde{h}\in \tilde{W}\colon\,\vec{T}(h) = \vec{T}^*_\infty\}.
\end{equation}
consisting of all graphons that meet the constraint.

In order to characterise the asymptotic behavior of the two ensembles, the entropy function 
of a Bernoulli random variable is essential. For $u\in [0,1]$, let 
\begin{equation}
\label{Idef1}
I(u) := \tfrac{1}{2}u\log u +\tfrac{1}{2}(1-u)\log(1-u).
\end{equation}
Extend the domain of this function to the graphon space $W$ by defining  
\begin{equation}
\label{Idef}
I(h) := \int_{[0,1]^2} \dd x\, \dd y\,\,I(h(x,y))
\end{equation}
(with the convention that $0\log0:=0$). On the quotient space $(\tilde{W},\delta_{\square})$, define $I(\tilde{h}) 
= I(h)$, where $h$ is any element of the equivalence class $\tilde{h}$. (In order to keep the notation minimal, 
we use $I(\cdot)$ for both \eqref{Idef1} and \eqref{Idef}. Below it will always be clear which of the two is being considered.) The function $\tilde{h} \mapsto I(\tilde{h})$ is well defined on the quotient space $\tilde{W}$ 
(see \cite[Lemma 2.1]{CV11}). Moreover, $\tilde{h} \mapsto I(\tilde{h}) + \tfrac12\log 2$ is the rate function of 
the large deviation principle in \cite{CV11}. 

The key result in \cite{dHMRS18} is the following variational formula for $s_\infty$ ($\cdot$ denotes the inner 
product for vectors).

\begin{theorem} {\rm \cite{dHMRS18}}
\label{th:Limit}
Subject to \eqref{softconstraint2}, \eqref{eq: Assumption T} and \eqref{eq:Assumption},
\begin{equation}
\label{sinftydef}
\lim_{n\to\infty} n^{-2} S_n(\Pmic \mid \Pcan) =: s_\infty
\end{equation}
with
\begin{equation}
\label{varreprsinfty}
s_\infty = \sup_{\tilde{h}\in \tilde{W}} \left[\vec{\theta}^*_\infty\cdot\vec{T}(\tilde{h})-I(\tilde{h})\right]
-\sup_{\tilde{h}\in \tilde{W}^*} \left[\vec{\theta}^*_\infty\cdot\vec{T}(\tilde{h}) - I(\tilde{h})\right].
\end{equation}
\end{theorem}

\noindent
Theorem~\ref{th:Limit} and the compactness of $\tilde{W}^*$ give us a \emph{variational 
characterisation} of ensemble equivalence: $s_\infty = 0$ if and only if at least one of the 
maximisers of $\vec{\theta}^*_\infty\cdot\vec{T}(\tilde{h})-I(\tilde{h})$ in $\tilde{W}$ also lies 
in $\tilde{W}^* \subset \tilde{W}$, i.e., satisfies the constraint $T^*_\infty$. 


\subsection{Edges and triangles}
\label{Sedtr} 

Theorem~\ref{th:Limit} allows us to identify examples where ensemble equivalence holds 
($s_\infty=0$) or is broken ($s_\infty>0$). In \cite{dHMRS18} a detailed analysis was given 
for the special case where the constraint is on the densities of edges and triangles. 

\begin{figure}[htbp]
\centering
\vspace{0.5cm}\hspace{1.5cm}
\includegraphics[width=0.7\linewidth,height=7cm]{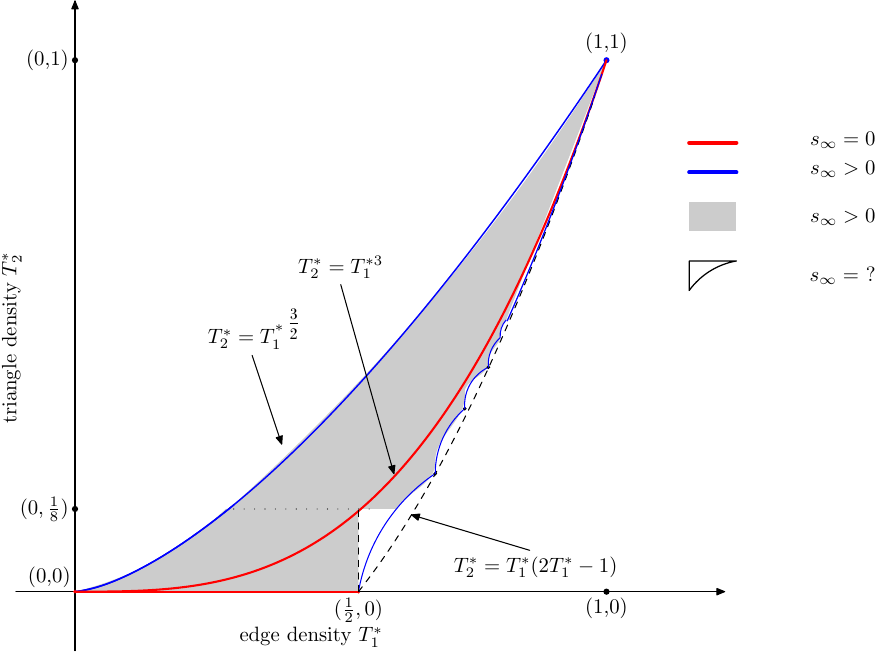}
\caption{\small The admissible edge-triangle density region is the region on and 
between the blue curves \cite{RS15}.}
\label{fig-scallopy}
\end{figure}

\begin{theorem} {\rm \cite{dHMRS18}}
\label{thm:equivalence}
For the edge-triangle model, $s_\infty=0$ when
\begin{itemize}
\item [$\circ$]
$T_2^* = T_1^{*{3}}$, 
\item [$\circ$]
$0<T_1^* \leq \frac{1}{2}$ and $T_2^* = 0$,
\end{itemize}
while $s_\infty>0$ when 
\begin{itemize}
\item [$\circ$]
$T_2^* \neq T_1^{*{3}}$ and $T_2^* \geq \tfrac18$,
\item [$\circ$]
$T_2^*\neq T_1^{*{3}}$, $0<T_1^*\leq\tfrac12$ and $0<T_2^*<\frac{1}{8}$,
\item [$\circ$]
$(T_1^*,T_2^*)$ lies on the scallopy curve $($the lower blue curve Figure~{\rm \ref{fig-scallopy}}$)$. 
\end{itemize}
\end{theorem} 

\noindent
Here, $T_1^*,T_2^*$ are in fact the limits $T_{1,\infty}^*,T_{2,\infty}^*$ in \eqref{eq: Assumption T}, but 
in order to keep the notation light we suppress the index $\infty$.

\noindent
Theorem~\ref{thm:equivalence} is illustrated in Figure~\ref{fig-scallopy}. The region on and between 
the blue curves, called the \emph{admissible region}, corresponds to the choices of $(T_1^*,T_2^*)$ 
that are \emph{graphical}, i.e., there exists a graph with edge density $T_1^*$ and triangle density 
$T_2^*$. The red curves represent ensemble equivalence, while the blue curves and the grey region 
represent breaking of ensemble equivalence. In the white region between the red curve and the lower 
blue curve we do not know whether there is breaking of ensemble equivalence or not (although we 
expect that there is). Breaking of ensemble equivalence arises from \emph{frustration} between the 
values of $T_1^*$ and $T_2^*$. In line with \cite{dHMRS18}, by frustration we mean that these two densities do not lie on the Erd\H{o}s-R\'enyi line. 

\begin{figure}[htbp]
\centering
\includegraphics[width=0.6\linewidth,height=6cm]{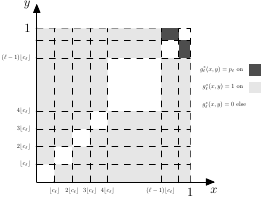}
\hspace{0.0cm}
\caption{\small $g_{\ell}^*$ for $\ell\in\mathbb{N}$ and $T_1^*\in(\frac{\ell-1}{\ell},\frac{\ell}{\ell+1}]$.}
\label{fig: graphon ell}
\end{figure}

The lower blue curve, called the \emph{scallopy curve}, consists of infinitely many pieces labelled
by $\ell \in \N \setminus \{1\}$. The $\ell$-th piece corresponds to $T_1^* \in (\frac{\ell-1}{\ell},
\frac{\ell}{\ell+1}]$ and a value $T_2^*$ that is a computable but non-trivial function of $T_1^*$, 
explained in detail in \cite{PR12,RS13,RS15}. The structure of the graphs drawn from the microcanonical 
ensemble was determined in \cite{PR12,RS15}. On the $\ell$-th piece:  
\begin{itemize}
\item
The vertex set can be partitioned into $\ell$ subsets. The first $\ell-1$ subsets have size $\lfloor c_{\ell}n\rfloor$, 
the last subset has size between $\lfloor c_\ell n\rfloor$ and $2\lfloor c_\ell n\rfloor$, with
\begin{equation}
\label{eq: for c}
c_{\ell} := \tfrac{1}{\ell+1}\left[1+ \sqrt{1-\tfrac{\ell+1}{\ell}\,T_1^*}\,\right] 
\in [\tfrac{1}{\ell+1},\tfrac{1}{\ell}).
\end{equation}
\item
The graph has the form of a complete $\ell$-partite graph, with some additional edges on the last subset that 
create no triangles within that last subset. 
\item
The optimal graphons have the form
\small
\begin{equation}
\label{eq: g*ell}
g^*_{\ell}(x,y) := \left\{
\begin{array}{ll}
1, &\exists\,1 \leq k < \ell\colon\,x~<kc_\ell<~y \mbox{ or }  y<kc_\ell<x,\\ 
p_\ell, &(\ell-1)c_\ell~<x<~\tfrac12[1+(\ell-1)c_\ell]<y \mbox{ or } 
(\ell-1)c_\ell~<y<~\tfrac12[1+(\ell-1)c]<x,\\
0, &\mbox{otherwise},
\end{array}
\right.
\end{equation}
\normalsize
where 
\begin{equation}
\label{eq: for p}
p_\ell = \frac{4c_\ell(1-\ell c_\ell)}{(1-(\ell-1)c_\ell)^2} \in (0,1].
\end{equation}
\end{itemize}
Figure~\ref{fig: cl pl T1} plots the coefficients $c_{\ell}$ and $p_{\ell}$ as a function of $T_1^*$ for $\ell\in\N$. 
Figure \ref{fig: graphon ell} plots the graphon $g^*_{\ell}$ for $\ell\in\N$. 

\begin{figure}[htbp]
\centering
\includegraphics[width=0.5\linewidth,height=6cm]{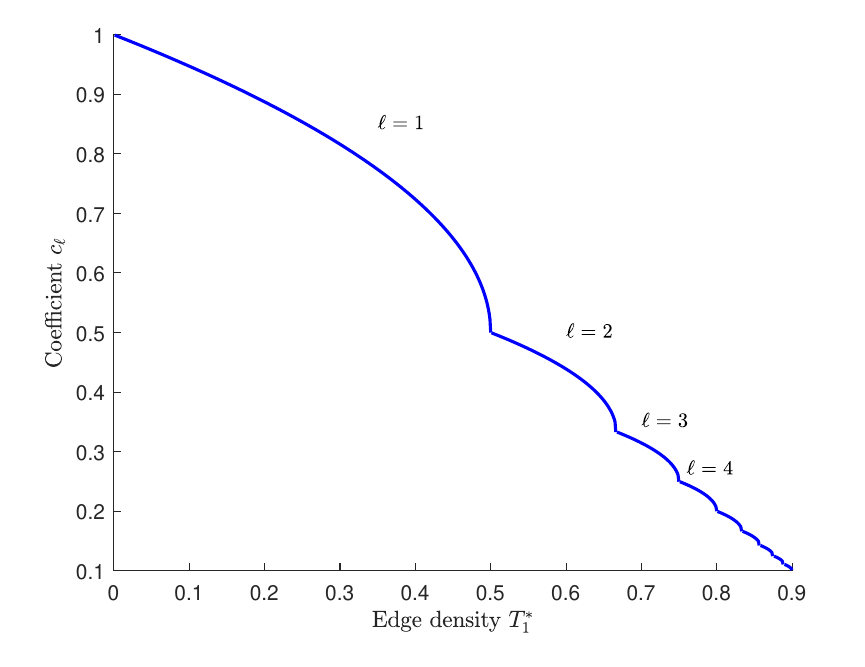}
\hspace{0.0cm}
\includegraphics[width=0.45\linewidth,height=6cm]{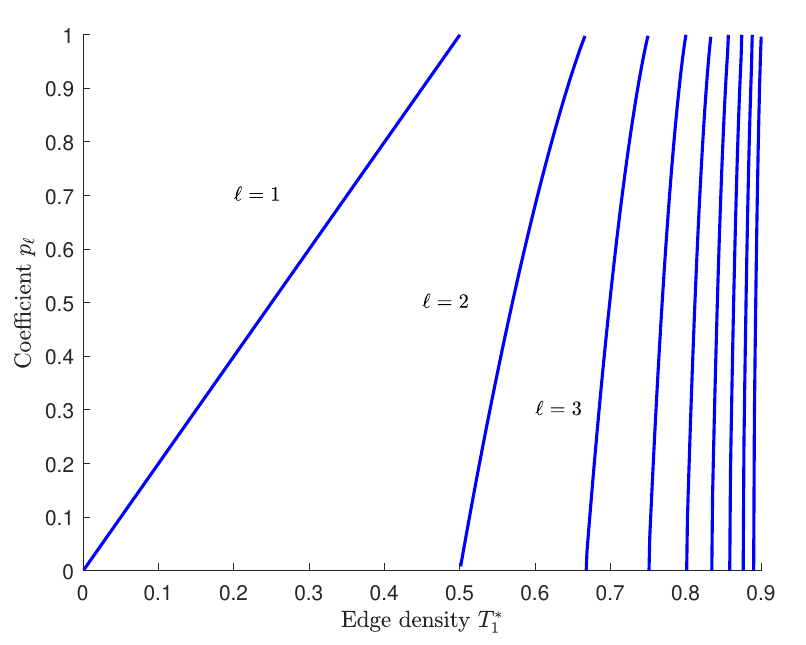}
\caption{\small $c_{\ell}$ (left) and $p_{\ell}$ (right) as a function of $T_1^*$.}
\label{fig: cl  pl T1}
\end{figure}

\section{Main results}
\label{S1.4}

In Section~\ref{ss.ass} we state two assumptions. In Section~\ref{ss.scal} we identify the scaling behaviour 
of $s_\infty$ for fixed $T_1^*$ and $T_2^* \downarrow T_1^{*3}$, respectively, $T_2^* \uparrow T_1^{*3}$. 
It turns out that the way in which $s_\infty$ tends to zero differs in these two cases. In Section~\ref{ss.opt} 
we characterise the asymptotic structure of random graphs drawn from the microcanonical ensemble for 
fixed $T_1^*$ and $T_2^* \downarrow T_1^{*3}$, respectively, $T_2^* \uparrow T_1^{*3}$. It turns out that 
the structure differs in these two cases.


\subsection{Assumptions}
\label{ss.ass}

Throughout the sequel we make the following two assumptions:

\begin{assumption}
\label{assumption}
{\rm Fix the edge density $T_1^*\in(0,1)$ and consider the triangle density $T_1^{*3}+\e$ for some $\e$, 
either positive or negative. Consider the associated Lagrange multipliers $\vec{\theta}^*_{\infty}
(\e):=(\theta^*_1(\e),\theta^*_2(\e))$. Then, for $\e$ sufficiently small, we have the representation 
\begin{equation}
\sup_{\tilde{h}\in\tilde{W}}\left[\theta_1^*(\e)T_1(\tilde{h}) + \theta_2^*(\e)T_2(\tilde{h}) - I(\tilde{h})\right] 
= [\theta_1T_1^* - I(T_1^*)] + (\gamma_1T_1^* + \gamma_2T_1^{*3})\e +O(\e^2),
\end{equation}
where $\theta_1 :=\theta^*_1(0), \gamma_1 = \theta^{*'}_1(0)$ and $\gamma_2 = \theta^{*'}_2(0)$. (It 
turns out that $\theta_2 :=\theta^*_2(0)=0$.)}
\end{assumption}

\begin{assumption}
\label{assumption 2}
{\rm Fix the edge density $T_1^*\hspace{-1pt}\in\hspace{-1pt}(0,1)$ and consider the triangle density 
$T_1^{*3}+\e$ for some $\e$, either positive or negative. Consider the microcanonical entropy 
\begin{equation}
\label{eq: micro ass}
-J(\e) := \sup\{-I(\tilde{h}):~\tilde{h}\in\tilde{W},~T_1(\tilde{h}) = T_1^*, ~T_2(\tilde{h}) = T_1^{*3}+\e\}.
\end{equation}
Then, for $\e$ sufficiently small, the maximizer $h_{\e}^*$ of \eqref{eq: micro ass} has the form 
\begin{equation}
\label{eq: opt graphon}
h_{\e}^* = T_1^* + g_{\e}, \qquad g_{\e} = g_{11}1_{I\times I} + g_{12}1_{(I\times J)\cup( J\times I)} 
+ g_{22}1_{J\times J},
\end{equation}
for some $g_{11},g_{12},g_{22}\in[-T_1^*,1-T_1^*]$ and $I,J\subset[0,1]$. }
\end{assumption}

\begin{remark}
{\rm Assumption 1 is needed to prove Theorems~\ref{thm:perturbationdown}--\ref{thm: perturbationup scallop}.
In Section \ref{S2.1 Chapter 5} we show that Assumption~\ref{assumption} is true when $T_1^*\in
[\frac{1}{2},1)$, which implies the scalings in \eqref{eq: el 0} and \eqref{eq: el 2} below. We can prove the 
scalings in \eqref{eq: el 0} and \eqref{eq: el 1} below without Assumption~\ref{assumption}, but at the cost 
of replacing `$=$' by `$\geq$' and replacing $\lim$ by $\limsup$. If Assumption \ref{assumption} fails, then 
these scalings hold with strict inequality.}
\end{remark}

\begin{remark}
{\rm Assumption 2 is needed to prove Propositions~\ref{thm:optimal perturbation2}--\ref{thm:optimal perturbation3}. Importantly, the validity of Assumption \ref{assumption 2} is firmly backed by the extensive numerical experiments performed 
in \cite{KRRS172}, showing that close to the ER line the optimizing graphon has the form  
\eqref{eq: opt graphon}.
The assumption reflects the following informal argument. Suppose that we want to 
maximise the microcanonical entropy among block graphons. Then we expect the entropy to decrease 
when we add more structure to the graphon. An $m \times m$ block graphon corresponds to a random 
graph where the vertices are divided into $m$ groups, and {\em within} each group we have an 
Erd\H{o}s-R\'enyi random graph with a certain retention probability. We expect that the microcanonical 
entropy decreases as $m$ increases. } 
\end{remark}


\subsection{Scaling of the specific relative entropy}
\label{ss.scal}

The variational problem $J(\e)$ in \eqref{eq: micro ass} was solved in \cite{KRRS18} for the case 
$\epsilon>0$, while the case $\epsilon$ remained unsolved. We consider small $\epsilon$ only, 
which is simpler and yields more intuition about the way the constraint is achieved. In \cite{KRRS17} 
the maximisers of the microcanonical entropy are identified numerically. The optimal graphons obtained 
agree with the optimal graphons that we derive rigorously. 

\begin{theorem}
\label{thm:perturbationdown}
For $T_1^* \in (0,1)$ and $T_1^*\neq \frac{1}{2}$,
\begin{equation}
\label{eq: el 0}
\lim_{\epsilon \downarrow 0}\, \epsilon^{-1} s_\infty(T_1^*,T_1^{*3} +3T_1^*\epsilon) 
= A(T_1^*) := - \frac{1}{1-2T_1^*}\log\frac{T_1^*}{1-T_1^*}  \in (0,\infty).
\end{equation}
\end{theorem}

\begin{theorem}
\label{thm:perturbationup}
For $T_1^* \in (0,\frac{1}{2}]$, 
\begin{equation}
\label{eq: el 1}
\lim_{\epsilon \downarrow 0}\, \epsilon^{-2/3} s_\infty(T_1^*,T_1^{*3} - T_1^{*3}\epsilon) 
= B(T^*_1) := \frac{1}{4}\frac{T_1^*}{1-T_1^*} \in (0,\infty). 
\end{equation}
\end{theorem}

\begin{theorem}
\label{thm: perturbationup scallop}
For $T_1^*\in(\frac{1}{2},1)$, 
\begin{equation}
\label{eq: el 2}
\lim_{\e\downarrow0}\e^{-2/3} s_{\infty}(T_1^*,T_1^{*3}-T_1^{*3}\e) = f(T_1^*,y^*)\in(0,\infty),
\end{equation}
where $y^* = y^*(T^*_1) \in (-T_1^*,0)$ is the unique point where the function $x\mapsto f(T_1^*,x)$ defined 
by 
\begin{equation}
\label{def: f function}
f(T_1^*,x):= T_1^{*2}\frac{I(T_1^*+x)-I(T_1^*) - I'(T_1^*)x}{x^2}, \quad x\in(-T^*_1,0),
\end{equation}
attains its global minimum.
\end{theorem}

The main implication of Theorems~\ref{thm:perturbationdown}--\ref{thm: perturbationup scallop} is
that, for a fixed value of the edge density, it is less costly in terms of relative entropy to increase the 
density of triangles than to decrease the density of triangles.  Above the ER line the cost is linear in 
the distance, below the ER line the cost is proportional to the $\tfrac23$-power of the distance.

\begin{figure}[htbp]
\centering
\includegraphics[width=0.46\linewidth,height=7cm]{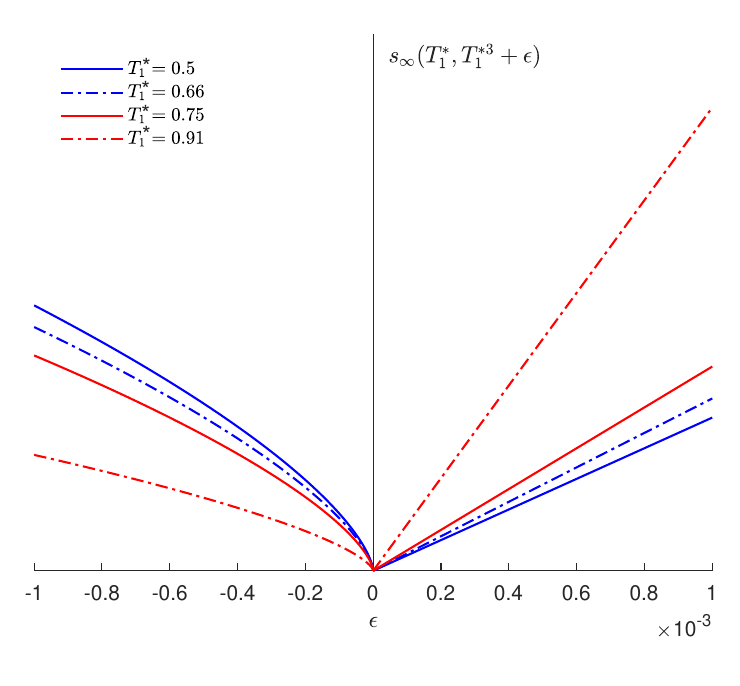}
\caption{\small Limit of scaled $s_{\infty}$ as a function of $\e$ for $\e$ sufficiently small.}
\label{fig: scaling entropy}
\end{figure}

For the case $T_1^* \in [\frac{1}{2},1)$ the above results are illustrated in Figure \ref{fig: scaling entropy}. 
In the left panel we plot $\epsilon^{2/3} f(T_1^*,y^*)$, in the right panel we plot $\epsilon A(T^*_1)$. In both 
panels we take $\e$ sufficiently small and pick four different values of $T_1^*$.

\subsection{Scaling of the optimal graphon}
\label{ss.opt} 

In Propositions~\ref{thm:optimal perturbation1}--\ref{thm:optimal perturbation3} below
we identify the structure of the optimal graphons corresponding to the perturbed constraints 
in the microcanonical ensemble in the limit as $n\to\infty$.

\begin{proposition}
\label{thm:optimal perturbation1}
For the pair of constraints $(T_1^*,T_1^{*3}+3T_1^*\e)$ with $\e>0$ sufficiently small, the 
unique optimiser $h_{\e}^*$ of the microcanonical entropy is given by
\begin{equation}
\label{eq: optimal g*}
h_{\e}^*(x,y) = \left\{ \begin{array}{ll}
\hspace{2pt}h_{11} +o(\epsilon), &(x,y) \in [0,\lambda\e]^2 ,\\[0.2cm]
1-T_1^*+h_1\e + o(\epsilon), &(x,y) \in [0,\lambda\e]\times (\lambda\e,1] \cup (\lambda\e,1]\times [0,\lambda\e],\\[0.2cm]
T_1^*+h_2\e + o(\epsilon), &(x,y) \in (\lambda\e,1]^2,\\
\end{array} \right.
\end{equation}
where 
\begin{equation}
\lambda :=\frac{1}{(1-2T_1^*)^2}, \qquad h_{1}:=\frac{1}{2}h_{2}, \qquad \qquad h_{2}:=-\frac{2}{1-2T_1^*},
\end{equation}
and $h_{11}$ solves the equation $I'(h_{11})=3I'(1-T_1^*)$. 
\end{proposition}

\begin{proposition}
\label{thm:optimal perturbation2}
For $T_1^*\in(0,\frac{1}{2}]$ and $T_2^* = T_1^{*3}-3T_1^*\e$, $\e \downarrow 0$, the optimal graphon is given by
\begin{equation}
\label{hpert2*}
h^*_{\epsilon} = T_1^* + \e^{1/3}{g}^* + o(\e^{2/3}) 
\qquad \text{\em (global perturbation)}
\end{equation}   
with $g^*$ defined by
\begin{equation}
\label{eq: optimal pert g* 2}
g^*(x,y) = \left\{ \begin{array}{ll}
\hspace{-2pt}-T_1^*, &(x,y) \in [0,\frac{1}{2}]^2 ,\\[0.2cm]
\hspace{5pt}T_1^*, &(x,y) \in [0,\frac{1}{2}]\times (\frac{1}{2},1] \cup (\frac{1}{2},1]\times [0,\frac{1}{2}],\\[0.2cm]
\hspace{-2pt}-T_1^*, &(x,y) \in (\frac{1}{2},1]^2.\\
\end{array} \right.
\end{equation}
\end{proposition}

\begin{proposition}
\label{thm:optimal perturbation3}
For $T_1^*\in(\frac{1}{2},1)$ and $T_2^* = T_1^{*3}-3T_1^*\e$, $\e\downarrow 0$, the optimal graphon is given by
\begin{equation}
\label{hpert2* scallop}
h^*_{\epsilon}  = T_1^*+g_{\e}^*
\qquad \text{\em (local perturbation)}
\end{equation}   
with $g_{\e}^*$ defined by 
\begin{equation}
g_{\e}^*(x,y) :=
\begin{cases}
~\frac{T_1^{*2}}{y*}\e^{2/3},~&\qquad (x,y)\in[0,1-\frac{T_1^*}{|y^*|}\e^{1/3}]^2,\\[0.2cm] 
~T_1^*\e^{1/3}, ~&\qquad (x,y)\in [0,1-\frac{T_1^*}{|y^*|}\e^{1/3}]\times[1-\frac{T_1^*}{|y^*|}\e^{1/3},1],\\[0.2cm]
&\qquad \text{ {\rm or} } (x,y)\in[1-\frac{T_1^*}{|y^*|}\e^{1/3},1]\times[0,1-\frac{T_1^*}{|y^*|}\e^{1/3}],\\[0.2cm]
~y^*,~&\qquad (x,y)\in[1-\frac{T_1^*}{|y^*|}\e^{1/3},1]^2,
\end{cases}
\end{equation} 
where $y^*\in(-T_1^*,0)$ is as in Theorem~\ref{thm: perturbationup scallop}.
\end{proposition}

The main implication of Propositions~\ref{thm:optimal perturbation1}--\ref{thm:optimal perturbation3} 
is that the optimal perturbation of the Erd\H{o}s-R\'enyi graphon is \emph{global} above the ER line 
and also below the ER line when the edge density is less than or equal to $\tfrac{1}{2}$, whereas it 
is \emph{local} below the ER line when the edge density is larger than $\tfrac{1}{2}$.


\section{Proofs of Theorems \ref{thm:perturbationdown}--\ref{thm: perturbationup scallop}}
\label{S2 Chapter 5}

In Sections \ref{S2.1 Chapter 5}--\ref{S2.3} we prove Theorems 
\ref{thm:perturbationdown}--\ref{thm: perturbationup scallop}, respectively. Along the way we use 
Propositions \ref{thm:optimal perturbation1}--\ref{thm:optimal perturbation3}, which we prove in 
Section \ref{Optimal per}.


\subsection{Proof of Theorem~\ref{thm:perturbationdown}}
\label{S2.1 Chapter 5}

Let
\begin{equation}
T_1(\e) = T_1^*, \qquad T_2(\e) = T_1^{*3} + 3T_1^{*3}\e.
\end{equation}
The factor $3T_1^*$ appearing in front of the $\e$ is put in for convenience. We know that for 
every pair of graphical constraints $(T_1(\e),T_2(\e))$ there exists a unique pair of Lagrange 
multipliers $(\theta_1(\e),\theta_2(\e))$ corresponding to these constraints (to ease the notation 
we suppress $*$). For an elaborate discussion on existence and uniqueness we refer the reader 
to \cite{dHMRS18}. By considering the Taylor expansion of the Lagrange multipliers $(\theta_1(\e),
\theta_2(\e))$ around $\e=0$, we obtain
\begin{equation}
\theta_1(\epsilon) = \theta_1 + \gamma_1 \epsilon + \tfrac{1}{2}\Gamma_1\e^2+O(\epsilon^3), 
\qquad 
\theta_2(\epsilon) = \gamma_2\epsilon+\tfrac{1}{2}\Gamma_2\e^2 + O(\epsilon^3),
\end{equation}
where
\begin{equation}
\theta_1(0)=\theta_1 = I'(T_1^*), \:~ \gamma_1 = \theta'_1(0), \:~ \Gamma_1 = \theta''_1(0), 
\:~ \theta_2(0)=0, \:~ \gamma_2 = \theta'_2(0), \:~ \Gamma_2 = \theta''_2(0).
\end{equation}
These equalities follow from \cite[Lemma 5.3]{dHMRS18}. For $\epsilon = 0$ we have $T_2^*(0) = T_1^{*3}$, 
which shows that the constraints correspond to those of the Erd\H{o}s-R\'enyi random graph. We denote the 
two terms in the expression for $s_\infty$ in \eqref{varreprsinfty} by $I_1,I_2$, 
i.e., $s_\infty = I_1-I_2$ with
\begin{equation}
I_1:=\sup_{\tilde{h}\in \tilde{W}} \big[\vec{\theta}_\infty\cdot\vec{T}(\tilde{h})-I(\tilde{h})\big],
\:\:\:I_2:=\sup_{\tilde{h}\in \tilde{W}^*} \big[\vec{\theta}_\infty\cdot\vec{T}(\tilde{h}) - I(\tilde{h})\big],
\end{equation}
and let $s_{\infty}(\epsilon)$ denote the relative entropy corresponding to the perturbed constraints. 
We distinguish between the cases $T_1^*\in[\frac{1}{2},1)$ and $T_1^*\in(0,\frac{1}{2})$.

\paragraph{Case I $T_1^*\in[\frac{1}{2},1)$:} 
According to \cite[Section 5]{dHMRS18}, if $T_1^*\in[\frac{1}{2},1)$ and $T_2^*\in[\frac{1}{8},1)$, then 
the corresponding Lagrange multipliers $(\theta_1,\theta_2)$ are both non-negative. Hence by 
\cite[Theorem 4.1]{CD13} we have that 
\begin{equation}
I_1 :=\sup_{\tilde{h}\in\tilde{W}}\left[\theta_1(\e)T_1(\tilde{h}) + \theta_2(\e) T_2(\tilde{h}) 
- I(\tilde{h})\right] = \sup_{0\leq u\leq 1}\left[\theta_1(\e)u+\theta_2(\e)u^3-I(u)\right],
\end{equation} 
and, consequently,
\begin{equation}
\label{eq: I1 CAseI}
I_1 = \theta_1(\e)u^*(\e) + \theta_2(\e)u^*(\e)^3-I(u^*(\e)).
\end{equation}
The optimiser $u^*(\e) \in (0,1)$ corresponding to the perturbed multipliers $\theta_1(\e)$ and 
$\theta_2(\e)$ is analytic in $\e$, as shown in \cite{RY13}. Therefore, a Taylor 
expansion around $\e=0$ gives 
\begin{equation}
\label{eq: expansion}
u^*(\e) = T_1^* + \delta \e + \tfrac12\Delta \e^2 +O(\e^3),
\end{equation}
where $\delta = {u^{*}}'(0)$ and $\Delta = {u^{*}}''(0)$. Hence $I_1$ can be written as 
\begin{equation}
\label{eq: expansion I}
I_1 = \theta_1T_1^* -I(T_1^*)+(\gamma_1 T_1^*+\gamma_2 T_1^{*3})\e + O(\e^2).
\end{equation}
Moreover,
\begin{equation}
\begin{aligned}
I_2 &= \big[\theta_1 + \gamma_1 \e + \tfrac{1}{2}\Gamma_1\e^2 +O(\e^3)\big]T_1^* 
+ \big[\gamma_2\e +\tfrac{1}{2}\Gamma_2\e^2 + O(\e^3)\big](T_1^{*3} + 3T_1^*\e) 
- \inf_{\tilde{h}\in\tilde{W}_{\e}^*} I(\tilde{h}) \\
&=\theta_1 T_1^* + \gamma_1 T_1^*\e +\tfrac{1}{2}\Gamma_1 T_1^*\e^2 
+ T_1^{*3} \gamma_2\e + \tfrac{1}{2}\Gamma_2 T_1^{*3}\e^2 
+3T_1^* \gamma_2 \e^2 - J^{\downarrow}(\e) + O(\e^3),
\end{aligned}
\end{equation}
where 
\begin{equation}
J^{\downarrow}(\e) :=\inf_{\tilde{h}\in\tilde{W}_{\e}^*} I(\tilde{h}), \qquad 
\tilde{W}_{\e}^*:=\{\tilde{h}\in \tilde{W}\colon\, T_1(\tilde{h})=T_1^*,\ T_2(\tilde{h}) = T_1^{*3} +3T_1^*\e\}.
\end{equation}
Consequently,
\begin{equation}
\label{eq: canonical entropy inf}
s_{\infty}(T_1^*,T_1^{*3}+3T_1^{*3}\e) = J^{\downarrow}(\e) - I(T_1^*)+O(\e^2).
\end{equation}
A straightforward computation of the entropy of $h_{\e}^*$ in \eqref{eq: optimal g*} shows that
\begin{equation}
\begin{aligned}
J^{\downarrow}(\e) = I(T_1^*)+I'(T_1^*)h_2\epsilon+o(\epsilon) 
= I(T_1^*) - \frac{1}{1-2T_1^*}\log\frac{T_1^*}{1-T_1^*}\e + o(\e).
\end{aligned}
\end{equation}
Hence we obtain that
\begin{equation}
s_{\infty}(T_1^*,T_1^{*3}+3T_1^{*3}\e) = - \frac{1}{1-2T_1^*}\log\frac{T_1^*}{1-T_1^*}\,\e  +o(\e).
\end{equation}

\paragraph{Case II $T_1^* \in(0,\frac{1}{2})$:}
Consider the term 
\begin{equation}
I_1 :=\sup_{\tilde{h}\in\tilde{W}}\left[\theta_1(\e)T_1(\tilde{h}) + \theta_2(\e) 
T_2(\tilde{h}) - I(\tilde{h})\right],
\end{equation}
as above. If Assumption \ref{assumption} applies, then this case can be treated in the same way as 
Case I. Otherwise, consider the lower bound
\begin{equation}
\sup_{\tilde{h}\in\tilde{W}}\left[\theta_1(\e)T_1(\tilde{h}) + \theta_2(\e) T_2(\tilde{h}) 
- I(\tilde{h})\right]\geq \sup_{0\leq u\leq 1}\left[\theta_1(\e)u+\theta_2(\e)u^3-I(u)\right].
\end{equation}
The arguments for Case I following \eqref{eq: I1 CAseI} still apply, and \eqref{eq: canonical entropy inf} 
is obtained with an inequality instead of an equality. 


\subsection{Proof of Theorem~\ref{thm:perturbationup}}
\label{S2.2}

In this section we omit the computations that are similar to those in the proof of 
Theorem~\ref{thm:perturbationdown}, as provided  in Section~\ref{S2.1 Chapter 5}. Let
\begin{equation}
T_1(\e) = T_1^*, \qquad T_2(\e) =  T_1^{*3} - T_1^{*3}\e.
\end{equation}
The factor $T_1^{*3}$ appearing in front of $\e$ is put in for convenience (without loss of generality). 
The perturbed Lagrange multipliers are
\begin{equation}
\theta_1(\epsilon) = \theta_1 + \gamma_1 \epsilon + \tfrac{1}{2}\Gamma_1\e^2+O(\epsilon^3), 
\qquad 
\theta_2(\epsilon) = \gamma_2\epsilon+\tfrac{1}{2}\Gamma_2\e^2 + O(\epsilon^3),
\end{equation}
where
\begin{equation}
\theta_1 = I'(T_1^*), \qquad \gamma_1 = \theta'_1(0), \qquad \Gamma_1 = \theta''_1(0) 
\qquad \gamma_2 = \theta'_2(0), \qquad \Gamma_2 = \theta''_2(0).
\end{equation}
We denote the two terms in the expression for $s_\infty$ in \eqref{varreprsinfty} by $I_1,I_2$, 
i.e., $s_{\infty} = I_1-I_2$, and let $s_{\infty}(\epsilon)$ denote the perturbed relative entropy. 
The computations for $I_1$ are similar as before, because the exact form of the constraint 
does not affect the expansions in \eqref{eq: expansion} and \eqref{eq: expansion I}. For $I_2$, 
on the other hand, we have 
\begin{equation}
\begin{aligned}
I_2 &= \theta_1T_1^*+\gamma_1T_1^*\e + \tfrac{1}{2}\Gamma_1 T_1^*\e^2
+T_1^{*3}\gamma_2\e+\tfrac{1}{2}\Gamma_2 T_1^{*3}\e^2-T_1^{*3}\gamma_2 \e^2 
- J_1^{\uparrow}(\e) +\,O(\epsilon^3) \\
&= \theta_1T_1^*+\gamma_1T_1^*\e +T_1^{*3}\gamma_2\e - J_1^{\uparrow}(\e) +O(\e^2),
\end{aligned}
\end{equation}
where 
\begin{equation}\label{eq:Jar}
J_1^{\uparrow}(\e) :=\inf_{\tilde{h}\in\tilde{W}_{\e}^*} I(\tilde{h}), \qquad 
\tilde{W}_{\e}^* := \big\{\tilde{h}\in \tilde{W}\colon\, T_1(\tilde{h})=T_1^*,\ T_2(\tilde{h}) = T_1^{*3} -T_1^{*3}~\e\big\}.
\end{equation}
Consequently,
\begin{equation}
\label{eq: entropy upwards 22}
s_{\infty}(T_1^*,T_1^*-T_1^{*3}\e) = J_1^{\uparrow}(\e) - I(T_1^*)+O(\e^2).
\end{equation}
Denote by $\tilde{h}_{\e}^*$ the optimiser of the variational problem 
$J_1^{\uparrow}(\e)$, as defined in \eqref{eq:Jar}. From Proposition \ref{thm:optimal perturbation2} we know that, for $T_1^*\in
(0,\frac{1}{2}]$, any optimal graphon in the equivalence 
class $\tilde{h}_{\e}^*$, denoted by $h_{\e}^*$ for simplicity in the notation, has the form
\begin{equation}
\label{hpert2* b}
h_{\e}^* = T_1^* + \e^{1/3}{g}^* + o(\e^{2/3}) 
\end{equation}   
with $g^*$ given by
\begin{equation}
\label{eq: optimal pert g* 2 b}
g^*(x,y) = \left\{ \begin{array}{ll}
\hspace{-2pt}-T_1^*, &(x,y) \in [0,\frac{1}{2}]^2 ,\\[0.2cm]
\hspace{5pt}T_1^*, &(x,y) \in [0,\frac{1}{2}]\times (\frac{1}{2},1] \cup (\frac{1}{2},1]\times [0,\frac{1}{2}],\\[0.2cm]
\hspace{-2pt}-T_1^*, &(x,y) \in (\frac{1}{2},1]^2.\\
\end{array} \right.
\end{equation}
Hence
\begin{equation}
J_1^{\uparrow}(\e) = I(T_1^*) +\frac{1}{2}T_1^{*2}I''(T_1^*)\e^{2/3} + o(\epsilon^{2/3}) 
= I(T_1^*) + \frac{1}{4}\frac{T_1^*}{1-T_1^*}\e^{2/3} + o(\epsilon^{2/3}),
\end{equation}
which gives
\begin{equation}
s_{\infty}(T_1^*,T_1^*-T_1^{*3}\e) = \frac{1}{4}\frac{T_1^*}{1-T_1^*} \e^{2/3} + o(\e^{2/3}).
\end{equation}


\subsection{Proof of Theorem~\ref{thm: perturbationup scallop}}
\label{S2.3}

The computations leading to the expression for the relative entropy in the right-hand side of 
\eqref{eq: el 2} are similar to those in Section \ref{S2.2}, and we omit them. Hence we have 
\begin{equation}
\label{eq: entropy upwards 2ell}
s_{\infty}(T_1^*,T_1^*-T_1^{*3}\e)= J_{2}^{\uparrow}(\e) - I(T_1^*)+O(\e^2),
\end{equation}
where, for $T_1^*\in(\frac{1}{2},1)$,
\begin{equation}\label{eq:Jar2}
J_{2}^{\uparrow}(\e) :=\inf_{\tilde{h}\in\tilde{W}_{\e}^*} I(\tilde{h}), \qquad 
\tilde{W}_{\e}^* := \big\{\tilde{h}\in \tilde{W}: T_1(\tilde{h})=T_1^*,\ T_2(\tilde{h}) = T_1^{*3} -T_1^{*3}~\e\big\}.
\end{equation}
Denote by $\tilde{h}_{\e}^*$ the optimiser of the variational problem 
$J_2^{\uparrow}(\e)$, as defined in \eqref{eq:Jar2}. From Proposition \ref{thm:optimal perturbation3} we know that, for $T_1^*\in
(\frac{1}{2},1)$, any optimal graphon in the equivalence 
class $\tilde{h}_{\e}^*$, denoted by $h_{\e}^*$ for simplicity in the notation, has the form
\begin{equation}
h_{\e}^* = T_1^* + g_{\e}^*
\end{equation}   
with $g_{\e}^*$ given by
\begin{equation}
g_{\e}^*(x,y) :=
\begin{cases}
~{\displaystyle \tfrac{T_1^{*2}}{y^*}}\e^{2/3},
~&\qquad (x,y)\in[0,1-{\displaystyle \tfrac{T_1^*}{|y^*|}}\e^{1/3}]\times [0,1-{\displaystyle \tfrac{T_1^*}{|y^*|}}\e^{1/3}],\\ 
\vspace{-4mm}&\\
~T_1^*\e^{1/3}, ~&\qquad (x,y)\in [0,1-{\displaystyle \tfrac{T_1^*}{|y^*|}}\e^{1/3}]
\times[1-{\displaystyle \tfrac{T_1^*}{|y^*|}}\e^{1/3},1]\\ \vspace{-4mm}&\\
&\qquad ~ \text{ {\rm or} } (x,y)\in[1-{\displaystyle \tfrac{T_1^*}{|y^*|}}\e^{1/3},1]
\times[0,1-{\displaystyle \tfrac{T_1^*}{|y^*|}}\e^{1/3}],\\ \vspace{-4mm}&\\
~y^*,~&\qquad (x,y)\in[1-{\displaystyle \tfrac{T_1^*}{|y^*|}}\e^{1/3},1]
\times[1-{\displaystyle \tfrac{T_1^*}{|y^*|}}\e^{1/3},1].
\end{cases}
\end{equation} 
Hence we have 
\begin{equation}
s_{\infty}(T_1^*,T_1^*-T_1^{*3}\e) = f(T_1^*,y^*)\e^{2/3} + o(\e^{2/3}),
\end{equation}
where $y^*\in(-T_1^*,0)$ minimizes the function $x\mapsto f(T_1^*,x)$ 
defined by 
\begin{equation}
\label{fdef}
f(T_1^*,x):= T_1^{*2}\frac{I(T_1^*+x)-I(T_1^*) - I'(T_1^*)x}{x^2}, \quad x\in(-T_1^*,0).
\end{equation}
We proceed by showing that $x\mapsto f(T_1^*,x)$ has a unique minimizer $y^*\in(-T_1^*,0).$

\noindent
$\circ$~First, we show that $f(T_1^*,x) >0$ for every $T_1^*\in(0,1)$ and every $x\in(-T_1^*,0)$, 
or equivalently 
\begin{equation}
I(T_1^*+x) - I(T_1^*) - I'(T_1^*)x >0.
\end{equation}
From the mean-value theorem we have that, for any given $x\in(-T_1^*,0)$, there exists $\xi\in(T_1^*+x,T_1^*)$ such that $I'(T_1^*+x) 
- I(T_1^*) = I'(\xi)x$. Because $I'$ is an increasing function, this implies that \begin{equation}
f(T_1^*,x) = (I'(\xi)-I'(T_1^*))x >0,
\end{equation}
recalling that $x\in(-T_1^*,0)$ and $\xi\in(T_1^*+x,T_1^*)$. 

\noindent $\circ$~Second, we show that the function $x\mapsto f(T_1^*,x)$ attains a unique minimum at some point $y^*\in(-T_1^*, 0)$. A straightforward computation shows that the derivative 
of $f(T_1^*, \cdot)$ is equal to 
\begin{equation}
\label{eq: derivative f}
f'(T_1^*, x) = T_1^{*2}\frac{\left(I'(T_1^*+x)-I'(T_1^*)\right)x^2 - 2x\left(I(T_1^*+x)-I(T_1^*)-I'(T_1^*)x\right)}{x^4}.
\end{equation}
Substituting $x=0$ we observe that $f'(T_1^*, 0) = \frac{1}{6}(T_1^*)^2\,I^{(4)}(T_1^*) >0$, and by taking the limit $x\rightarrow -T_1^*$ we observe that 
\begin{equation}
\lim_{x\rightarrow-T_1^*}f'(T_1^*, x) = -\infty.
\end{equation} 
Hence the function $f(T_1^*, \cdot)$ is decreasing in a neighborhood of $-T_1^*$, while it is increasing in a neighborhood of zero. Consequently, there must be at least one point in $(-T_1^*,0)$ where the derivative is zero.

\noindent $\circ$~It remains to show that there is a unique point $y^*$ in $(-T_1^*,0)$ where the derivative is zero. Suppose that $y^*$ is such a point where $f'(T_1^*, y^*) = 0$. Then from \eqref{eq: derivative f} we know that 
\begin{equation}
\label{eq: derivative f2}
\left(I'(T_1^*+y^*)-I'(T_1^*)\right)y^* - 2\left(I(T_1^*+y^*)-I(T_1^*)-I'(T_1^*)y^*\right) = 0.
\end{equation}
From the mean-value theorem we know that there exist $\xi_1\in(T_1^*+y^*, T_1^*)$ and $\xi_2\in(T_1^*+y^*, T_1^*)$ such that $I''(\xi_1)y^* = I'(T_1^*+y^*)-I'(T_1^*)$ and $I'(\xi_2)y^* = I(T_1^*+y^*)-I(T_1^*)$. Moreover, $\xi_2$ is unique since $I$ is a convex function. This is not necessarily true for $\xi_1$. Hence \eqref{eq: derivative f2} becomes
\begin{equation}
\label{eq: derivative f3}
I''(\xi_1)y^* - 2(I'(\xi_2)-I'(T_1^*)) = 0.
\end{equation}
Applying again the mean-value theorem, we get that there exists a, not necessarily unique, $\xi_3\in(\xi_2, T_1^*)$ such that $I'(\xi_2)-I'(T_1^*) = I''(\xi_3)(\xi_2-T_1^*)$. Substituting this into \eqref{eq: derivative f3}, we obtain
\begin{equation}
\label{eq: derivative f4}
I''(\xi_1)y^* - 2I''(\xi_3)(\xi_2 - T_1^*)= 0.
\end{equation}
Hence any possible solution $y^*$ satisfies the equation
\begin{equation}
\label{eq: derivative f5}
y^* =-\frac{2I''(\xi_3)(T_1^*-\xi_2)}{I''(\xi_1)} = -\frac{2\xi_1(1-\xi_1)(T_1^*-\xi_2)}{\xi_3(1-\xi_3)}.
\end{equation}
Multiple solutions may arise due to the non-uniqueness of $\xi_1$ and $\xi_3$. Notice that if the function $I'$ attains a given slope at some $\xi^*$, it does so as well at $1-\xi^*$ (and nowhere else); use that $I''(u)$ is given by $1/(u(1-u))$. Therefore, other possible solutions may occur when we replace $\xi_1$ by $1-\xi_1$ and/or $\xi_3$ by $1-\xi_3$. However, it is directly seen that the right-hand side of \eqref{eq: derivative f5} is invariant under this substitution. Hence the point $y^*$ where the derivative of $f(T_1^*, \cdot)$ is equal to zero is unique. We finalize the proof of Proposition \ref{thm:optimal perturbation3} in the next section.


\section{Proofs of Propositions \ref{thm:optimal perturbation1}--\ref{thm:optimal perturbation3}}
\label{Optimal per}

In Section \ref{S5.1} we state two lemmas (Lemmas \ref{lem: perturbation up optimal 1}--\ref{lem: perturbation 
up optimal 2} below) about variational formulas encountered in Section \ref{S2 Chapter 5}, and use them 
to prove Propositions \ref{thm:optimal perturbation1}--\ref{thm:optimal perturbation3}. In Section \ref{appA3} 
we provide the proof of these two lemmas, which requires an additional lemma about a certain function related
to \eqref{fdef} (Lemma~\ref{lemma: function T1} below), whose proof is deferred to Section~\ref{end}.


\subsection{Key lemmas}
\label{S5.1}

In Section \ref{S2 Chapter 5} the following variational problems were encountered:
\begin{itemize}
\item[(1)] For $T_1^*\in(0,1)$,
\begin{equation}
\label{eq: VP1}
J^{\downarrow}(\e) = \inf\big\{ I(\tilde{h})\colon\,\tilde{h}\in\tilde{W},\,
T_1(\tilde{h})= T_1^*,\,T_2(\tilde{h}) = T_1^{*3} + 3T_1^{*3}\e\big\}.
\end{equation}
\item[(2)] For $T_1^*\in (0,\frac{1}{2}]$,
\begin{equation}
\label{eq: VP2} 
J_1^{\uparrow}(\e) = \inf\big\{ I(\tilde{h})\colon\,\tilde{h}\in \tilde{W},\,
T_1(\tilde{h})= T_1^*,\,T_2(\tilde{h}) = T_1^{*3} -T_1^{*3}\e\big\}.
\end{equation}
\item[(3)] For $T_1^*\in (\frac{1}{2},1)$,
\begin{equation}
\label{eq: VP3} 
J_{2}^{\uparrow}(\e) = \inf\big\{ I(\tilde{h})\colon\,\tilde{h}\in \tilde{W},\,
T_1(\tilde{h}) = T_1^*,\,T_2(\tilde{h}) = T_1^{*3} -T_1^{*3}\e\big\}.
\end{equation}
\end{itemize}
The difference between the variational problems in \eqref{eq: VP2} and \eqref{eq: VP3} lies only in the possible 
values $T_1^*$ can take. We give them separate displays in order to easily refer to them later on. 

In order to prove Propositions \ref{thm:optimal perturbation1}--\ref{thm:optimal perturbation3}, we need to analyse 
the above three variational problems for $\e$ sufficiently small. The variational formula in \eqref{eq: VP1} was
already analysed in \cite{KRRS18}. Solving the equations given in \cite[Theorem 1.1]{KRRS18} for the case of triangle density equal to $T_1^{*3} + 3T_1^{*3}\epsilon$ and $\epsilon$ sufficiently small enough we obtain the graphon given in \eqref{eq: optimal g*}. In what follows we concentrate on the variational formulas in \eqref{eq: VP2} 
and \eqref{eq: VP3}. We analyse the latter with the help of a perturbation argument. In particular, we show that 
the optimal perturbations are those given in \eqref{hpert2*} and \eqref{hpert2* scallop}, respectively. The claims 
in Propositions \ref{thm:optimal perturbation2} and \ref{thm:optimal perturbation3} follow directly from the 
following two lemmas. 

We remind the reader that Assumption \ref{assumption 2} is in forse, i.e., we look for the optimal graphon in 
the class of two-step graphons. 

\begin{lemma}
\label{lem: perturbation up optimal 1}
Let $T_1^*\in(0,\frac{1}{2}]$. For $\e>0$ sufficiently small, 
\begin{equation}
\label{eq: VP2 solved} 
J_1^{\uparrow}(\e) = I(T_1^*) + \frac{1}{4}\frac{T_1^*}{1-T_1^*}\e^{2/3}+o(\e^{2/3}).
\end{equation}
\end{lemma}

\begin{lemma}
\label{lem: perturbation up optimal 2}
Let $T_1^*\in(\frac{1}{2},1)$. For $\e>0$ sufficiently small, 
\begin{equation}
\label{eq: VP3 solved} 
J_{2}^{\uparrow}(\e) = I(T_1^*) + f(T_1^*,y^*)\e^{2/3} + o(\e^{2/3}),
\end{equation}
where $f(T_1^*,x)$, $x\in(-T_1^*,0)$, and $y^*$ are as defined in Theorem {\rm \ref{thm: perturbationup scallop}}.
\end{lemma}

In what follows we use the notation $f(\e)\asymp g(\e)$ when ${f(\e)}/{g(\e)}$ converges to a positive finite 
constant as $\e\downarrow0$, and $f(\e)=\omega(g(\e))$ when ${f(\e)}/{g(\e)}$ diverges as $\e\downarrow0$.


\subsection{Proof of Lemmas \ref{lem: perturbation up optimal 1}--\ref{lem: perturbation up optimal 2}}
\label{appA3}

\begin{proof}
Instead of the variational problems ${J}_1^{\uparrow}(\e)$ and ${J}_2^{\uparrow}(\e)$ we will consider for 
$\e>0$ the variational problem
\begin{equation}
\label{eq: variational problem down}
{J}^{\uparrow}(\e) = \inf\{I(\tilde{h}): \tilde{h}\in\tilde{W}, T_1(\tilde{h}) = T_1^*, T_2(\tilde{h})=T_1^{*3}-T_1^{*3}\e\}.
\end{equation}
Below we provide the technical details leading to the optimal perturbation corresponding to \eqref{eq: variational 
problem down}. At some point we will distinguish between the two cases $T_1^*\in(0,\frac{1}{2})$ and 
$T_1^*\in(\frac{1}{2},1)$, yielding the optimisers for ${J}_1^{\uparrow}(\e)$ and ${J}_2^{\uparrow}(\e)$. We 
denote the optimiser of \eqref{eq: variational problem down} by $\tilde{h}_{\e}^{*\uparrow}$ (in order to keep 
the notation light, we denote a representative element by $h_{\e}^{*\uparrow}$). 

We start by writing the optimiser in the form $h_{\e}^{*\uparrow} = T_1^* + \Delta H_{\e}$ for some 
perturbation term $\Delta H_{\e}$, which must be a bounded symmetric function on the unit square $[0,1]^2$ 
taking values in $\mathbb{R}$. The optimiser $h_{\e}^{*\uparrow}$ has to meet the two constraints 
\begin{equation}
T_1(h_{\e}^{*\uparrow}) = T_1^*, \qquad \qquad T_2(h_{\e}^{*\uparrow}) 
= T_1^{*3} - T_1^{*3}\e.
\end{equation}
Consequently, $\Delta H_{\e}$ has to meet the two constraints
\begin{eqnarray}
\label{constraints DH1 D}
(K_1): \hspace{1.5cm} \mathcal{K}_1 &:=& \int_{[0,1]^2} \dd x\,\dd y ~\Delta H_{\e}(x,y) = 0,\\
\label{constraints DH2 D}
(K_2): \hspace{0.65cm} \mathcal{K}_{2}+\mathcal{K}_3 &:=& 
3T_1^*\int_{[0,1]^3} \dd x\,\dd y\,\dd z\, ~\Delta H_{\e}(x,y)\Delta H_{\e}(y,z)\\ \nonumber
&&+ \int_{[0,1]^3} \dd x\,\dd y\,\dd z\, ~\Delta H_{\e}(x,y)\Delta H_{\e}(y,z)
\Delta H_{\e}(z,x) =-T_1^{*3}\e. 
\end{eqnarray}
Observe that $\mathcal{K}_2 \geq 0$. By Assumption \ref{assumption 2}, we restrict to graphons of the 
form $T_1^*+\Delta H_{\e}$ with
\begin{equation}
\label{eq: perturbation DHe}
\Delta H_{\e} = g_{11}1_{I\times I} + g_{12}1_{(I\times J)\cup (J\times I)} + g_{22} 1_{J\times J},
\end{equation}
where $g_{11},g_{12},g_{22}\in(-T_1^*,1-T_1^*)$ and $I\subset[0,1]$, $J=I^c$. From \eqref{constraints DH1 D} 
we have
\begin{equation}
\label{eq: DDH1}
\mathcal{K}_1=\lambda(I)^2g_{11} + 2\lambda(I)(1-\lambda(I))g_{12} +(1-\lambda(I))^2g_{22} =0,
\end{equation}
which yields
\begin{equation}
\label{eq: g12 upp}
g_{12} = -\frac{1}{2}\left(\frac{\lambda(I)}{1-\lambda(I)}g_{11} + \frac{1-\lambda(I)}{\lambda(I)}g_{22}\right).
\end{equation}
A straighforward computation shows that 
\begin{equation}
\label{eq: second order integral}
\mathcal{K}_2=3T_1^*\left(\lambda(I)^3g_{11}^2+\left(1-\lambda(I)\right)^3g_{22}^2+2\lambda(I)(1-\lambda(I))g_{12}\left(\lambda(I)g_{11}+(1-\lambda(I))g_{22}+\frac{1}{2}g_{12}\right)\right).
\end{equation}
Using \eqref{eq: g12 upp}, we obtain
\begin{equation}
\label{eq: g12 upp 2}
\mathcal{K}_2=3T_1^*\frac{1}{4}\lambda(I)(1-\lambda(I))\left(\frac{\lambda(I)}{(1-\lambda(I))}g_{11} 
- \frac{1-\lambda(I)}{\lambda(I)}g_{22}\right)^2.
\end{equation}
With a similar reasoning we also obtain that 
\begin{equation}
\label{eq: third order integral}
\mathcal{K}_3=\lambda(I)^3g_{11}^3+(1-\lambda(I))^3g_{22}^3+3g_{12}^2\lambda(I)(1-\lambda(I))\Big(\lambda(I)g_{11}+g_{22}(1-\lambda(I))\Big).
\end{equation}
We claim that in order to find the optimal graphon corresponding to the microcanonical entropy, it suffices to 
solve the following equations:
\begin{equation}
\label{eq: three integral equations}
\mathcal{K}_1 =0, \qquad \mathcal{K}_2 = 0, \qquad \mathcal{K}_3 = -T_1^{*3}\e. 
\end{equation}
We prove this claim in Appendix~\ref{app}. The idea is that for $\e$ sufficiently small, if $\mathcal{K}_2>0$, 
then we cannot have $\mathcal{K}_2+\mathcal{K}_3<0$. From the second equation in \eqref{eq: three integral 
equations} we obtain
\begin{equation}
\label{eq: g11}
g_{11} = \frac{(1-\lambda(I))^2}{\lambda(I)^2}g_{22},
\end{equation}
and substitution into \eqref{eq: g12 upp} gives
\begin{equation}
\label{eq: g12}
g_{12} = -\frac{1-\lambda(I)}{\lambda(I)}g_{22}.
\end{equation}
The third equation in \eqref{eq: three integral equations} now yields 
\begin{equation}
\label{eq: triple}
g_{22}\frac{1-\lambda(I)}{\lambda(I)} =-T_1^* \e^{1/3}.
\end{equation}
Substituting this equation into \eqref{eq: g11} and \eqref{eq: g12}, we obtain the two relations
\begin{equation}
\label{eq: g11 and g12 2}
g_{11} = -\frac{1-\lambda(I)}{\lambda(I)}T_1^*\e^{1/3}, \qquad  g_{12} = T_1^* \e^{1/3}.
\end{equation}
In what follows we need to distinguish between the following three cases: 
\begin{itemize}
\item[\bf (I)] $\lambda(I)$ is constant and independent of $\epsilon$;
\item[\bf (II)] $\lambda(I) \asymp \e^{1/3}$;
\item[\bf (III)] $\lambda(I)=\omega(\e^{1/3})$.
\end{itemize}
The case $\lambda(I)=o(\e^{1/3})$ can be excluded, since it yields via \eqref{eq: g11 and g12 2} that 
$g_{11}$ diverges as $\e \downarrow 0$, while from \eqref{eq: perturbation DHe} we argued that 
$g_{11}\in(-T_1^*,1-T_1^*)$. Hence the above three cases are exhaustive. We treat them separately by 
computing their corresponding microcanonical entropies. Afterwards, by comparing the three entropies we 
identify the optimal graphon. A straightforward computation yields 
\begin{align*}
I(h_{\e}^*) &= \lambda(I)^2 I\left(T_1^*- T_1^*\frac{1-\lambda(I)}{\lambda(I)}\e^{1/3}\right) 
+ 2\lambda(I)(1-\lambda(I))I\left(T_1^*+T_1^*\e^{1/3}\right) \\
&+ (1-\lambda(I))^2I\left(T_1^* - T_1^*\frac{\lambda(I)}{1-\lambda(I)}\e^{1/3}\right). 
\numberthis{\label{eq:Taylor cases}}
\end{align*}

\paragraph{Case (I).}
A Taylor expansion in $\e$ of all three terms in \eqref{eq:Taylor cases} yields 
\begin{equation}
\label{eq: Micro 11}
I(h_{\e}^*) = I(T_1^*) + \frac{1}{2}I''(T_1^*)T_1^{*2}\e^{2/3} 
- \frac{1}{6}I^{(3)}(T_1^*)T_1^{*3}\left(\frac{(1-2\lambda(I))^2}{\lambda(I)(1-\lambda(I))}\right)\e + o(\e).
\end{equation}

\paragraph{Case (II).} 
We have $\lambda(I)=c\e^{1/3}+o(\e^{1/3})$ for some constant $c>0$. Substituting the expressions for 
$g_{11}$, $g_{12}$, $g_{22}$ from \eqref{eq: triple} and \eqref{eq: g11 and g12 2} into \eqref{eq:Taylor cases}, 
and looking only at the leading order term, we obtain 
\begin{equation}
\label{eq: Second entropy}
I(h_{\e}^*) = I(T_1^*) + \Bigg(c^2\left(I\left(T_1^*-\frac{T_1^*}{c}\right)-I(T_1^*)\right)
+cT_1^*I'(T_1^*)\Bigg)\e^{2/3}+ O(\e).    
\end{equation}

\paragraph{Case (III).} 
A Taylor expansion in $\e$ of all three terms in \eqref{eq:Taylor cases} yields
\begin{equation}
\label{eq: Third entropy}
I(h_{\e}^*) = I(T_1^*) +\frac{1}{2}T_1^{*2}I''(T_1^*)\omega(\e^{2/3}) + \omega(\e).
\end{equation}

In order to determine the optimal graphon we need to compare the expressions in 
\eqref{eq: Micro 11}--\eqref{eq: Third entropy}. The leading order term in \eqref{eq: Third entropy} 
is $\omega(\e^{2/3})$, which entails that for $\e$ sufficiently small this term is larger than the 
corresponding second order terms in \eqref{eq: Micro 11} and \eqref{eq: Second entropy}. Hence 
it suffices to compare the second order terms in \eqref{eq: Micro 11} and \eqref{eq: Second entropy}. 
For this we use the following lemma.

\begin{lemma}
\label{lemma: function T1}
Consider, for $T_1^*\in(0,1)$ fixed and $x>0$, the function 
\begin{equation}
f(T_1^*,x) := x^2\left(I\left(T_1^*-\frac{T_1^*}{x}\right)-I(T_1^*)\right)+xT_1^*I'(T_1^*).
\end{equation}
Then the following properties hold:
\begin{itemize}
\item[$\circ$]
If $T_1^*\in(0,\frac{1}{2})$, then $f(T_1^*,x)>\frac{1}{2}T_1^{*2}I''(T_1^*)$ for all $x>0$.
\item[$\circ$] 
If $T_1^*\in(\frac{1}{2},1)$, then there exists $x^*>0$ such that $f(T_1^*,x^*)<\frac{1}{2}
T_1^{*2}I''(T_1^*)$.
\end{itemize}
\end{lemma}

The proof is given in Section~\ref{end}. Using Lemma~\ref{lemma: function T1}, we find that, 
for $T_1^*\in(0,\frac{1}{2})$ and all $c>0$, the microcanonical entropy in \eqref{eq: Micro 11} 
is larger than the microcanonical entropy in \eqref{eq: Second entropy}, while for $T_1^*\in(\frac{1}{2},1)$ 
there exists $c^*>0$ such that the microcanonical entropy in \eqref{eq: Micro 11} is smaller than 
the microcanonical entropy in \eqref{eq: Second entropy}. To complete the proof we determine the 
optimal graphons. 

\paragraph{Optimal graphon for $T_1^*\in(0,\frac{1}{2})$.}

The microcanonical entropy computed in \eqref{eq: Micro 11} is equal to
\begin{equation}
I(h_{\e}^*) = I(T_1^*) + \frac{1}{2}I''(T_1^*)T_1^{*2}\e^{2/3} 
- \frac{1}{6}I^{(3)}(T_1^*)T_1^{*3}\left(\frac{(1-2\lambda(I))^2}{\lambda(I)(1-\lambda(I))}\right)\e + O(\e).
\end{equation}
Since $T_1^*\in(0,\frac{1}{2})$, we have $I^{(3)}(T_1^*)<0$. Hence, in order to find the optimal graphon 
we need to minimise 
\begin{equation}
\frac{(1-2\lambda(I))^2}{\lambda(I)(1-\lambda(I))}.
\end{equation}
The minimum is achieved for $\lambda(I)=\frac{1}{2}$, which yields the graphon in Proposition 
\ref{thm:optimal perturbation2}.

\paragraph{Optimal graphon for $T_1^*\in(\frac{1}{2},1)$.} 

The microcanonical entropy computed in \eqref{eq: Second entropy} is equal to 
\begin{equation}
I(h_{\e}^*) = I(T_1^*) + \Bigg(c^2\left(I\left(T_1^*-\frac{T_1^*}{c}\right)-I(T_1^*)\right)
+cT_1^*I'(T_1^*)\Bigg)\e^{2/3}+ O(\e).    
\end{equation}
Hence, in order to find the optimal graphon we need to minimise the second order term. This is done
in Lemma \ref{lemma: function T1} and yields the graphon in Proposition \ref{thm:optimal perturbation3}.
\end{proof}


\subsection{Proof of Lemma \ref{lemma: function T1}}
\label{end}

\begin{proof}
Via the substitution $y:={T_1^*}/{x}$ we see that we can equivalently work with the function, defined for 
$y\in(0,T_1^*)$,
\begin{equation}
\label{fff}
f\left(T_1^*,\frac{T_1^*}{y}\right)- \frac{1}{2}T_1^{*2}I''(T_1^*)
= T_1^{*2}\frac{I(T_1^*-y)-I(T_1^*)+yI'(T_1^*)-\frac{1}{2}y^2I''(T_1^*)}{y^2}, 
\end{equation}
which we write for simplicity as $(T_1^*/y)^2\,\check f(T_1^*,y)$, with $\check f(T_1^*,y)$ the numerator in 
right-hand side of \eqref{fff}. It suffices to prove that: (i)~$\check f(T_1^*,y)>0$ for $T_1^*\in(0,\frac{1}{2})$ 
and all $y>0$; (ii)~$\check f(T_1^*,y)<0$ for $T_1^*\in(\frac{1}{2},1)$ and some $y>0$.

Our next observation is that
\begin{equation}
\label{Tay}
\check f(T_1^*,y) =\sum_{k=0}^\infty I^{(k)}(T_1^*)\frac{(-y)^k}{k!}- I(T_1^*)
+yI'(T_1^*)-\frac{1}{2}y^2I''(T_1^*)=\sum_{k=3}^\infty I^{(k)}(T_1^*)\frac{(-y)^k}{k!}.
\end{equation}
An elementary computation shows that, for $t\in(0,1)$ and $k\in\N\setminus\{1\}$,
\begin{equation}
I^{(k)}(t) = \frac{(k-2)!}{2}\left(\frac{(-1)^k}{t^{k-1}}+\frac{1}{(1-t)^{k-1}}\right).
\end{equation}
For $k$ even, $I^{(k)}(t)>0$ for all $t\in(0,1)$. For $k$ odd, $I^{(k)}(t)<0$ for $t\in(0,\frac{1}{2})$ and $I^{(k)}(t)>0$ 
for $t\in(\frac{1}{2},1)$. The above properties imply that, for all $T_1^*\in(0,\frac{1}{2})$ and all $y>0$, 
$I^{(k)}(T_1^*)\,{(-y)^k}>0$, so that \eqref{Tay} immediately implies claim (i). Claim (ii) follows from 
\eqref{Tay} in combination with
\begin{equation}
\check f(T_1^*,0) = \check f^{(1)}(T_1^*,0)  = \check f^{(2)}(T_1^*,0)  =0,\:\:\: 
\check f^{(3)}(T_1^*,0)=-I^{(3)}(T_1^*)<0, 
\end{equation}
where $\check f^{(k)}(T_1^*,y)$ denotes the $k$-th derivative of $\check f(T_1^*,y)$ with respect to $y$.
\end{proof}


\appendix

\section{Appendix}
\label{app}

In this section we prove the claim made in \eqref{eq: three integral equations}. Let 
\begin{equation}
\begin{aligned}
\mathcal{K}_1 &= \lambda(I)^2g_{11}+2\lambda(I)(1-\lambda(I))g_{12}+(1-\lambda(I))^2g_{22},\\
\mathcal{K}_{2} &= 3T_1^*\left(\frac{1}{4}\lambda(I)(1-\lambda(I))\left(\frac{\lambda(I)}{1-\lambda(I)}g_{11}
-\frac{1-\lambda(I)}{\lambda(I)}g_{22}\right)^2\right),\\
\mathcal{K}_3 &= \lambda(I)^3g_{11}^3+(1-\lambda(I))^3g_{22}^3
+3g_{12}^2\lambda(I)(1-\lambda(I))\Big(\lambda(I)g_{11}+(1-\lambda(I))g_{22}\Big).
\end{aligned}
\end{equation}
From the constraints on the perturbation we know that 
\begin{equation}
\mathcal{K}_1 =0, \qquad \mathcal{K}_2+\mathcal{K}_3 = -T_1^{*3}\e.
\end{equation}
We will show that it suffices to solve $\mathcal{K}_2=0$ and $\mathcal{K}_3=-T_1^{*3}\e$. The argument 
we use is similar to the one used in Section \ref{appA3} to find the optimal graphon. Since $\mathcal{K}_2 + \mathcal{K}_3 = -T_1^{*3}\e$, and $\mathcal{K}_2\geq0$ independently of $\e$, it must be that $\mathcal{K}_3 
= -c\e$ for some constant $c>0$. Using \eqref{eq: g12 upp}, after some straightforward computations we 
obtain
\begin{equation}
\begin{aligned}
\mathcal{K}_3 &=\lambda(I)^3g_{11}^3+(1-\lambda(I))^3g_{22}^3
+\frac{3}{4}\frac{\lambda(I)^4}{1-\lambda(I)}g_{11}^3+\frac{3}{4}\frac{(1-\lambda(I))^4}{\lambda(I)}g_{22}^3\\
&\qquad +\frac{3}{2}\lambda(I)^2(1-\lambda(I))g_{11}^2g_{22}+\frac{3}{2}(1-\lambda(I))^2\lambda(I)g_{11}g_{22}^2
+\frac{3}{4}(1-\lambda(I))^3g_{11}g_{22}^2+\frac{3}{4}\lambda(I)^3g_{11}^2g_{22}.
\end{aligned}
\end{equation}
We need to consider the following four cases:
\begin{itemize}
\item[(1)] $\lambda(I)^3g_{11}^3\asymp -\e$.
\item[(2)] $\frac{1}{\lambda(I)}g_{22}^3\asymp-\e$.
\item[(3)] $\lambda(I)^2g_{11}^2g_{22}\asymp-\e$.
\item[(4)] $\lambda(I)g_{11}g_{22}^2\asymp-\e$.
\end{itemize}
These cases are exhaustive because they cover all possible ways a term in $\mathcal{K}_3$ can be asymptotically 
of the order $-\e$. We first observe that, because of symmetry in $\lambda(I),g_{11}$ and $g_{22}$, the term 
$(1-\lambda(I))^3g_{22}$ can be dealt with in a similar way as in case (1). We show that in all cases if 
$\mathcal{K}_2>0$, then $\mathcal{K}_2=\omega(\e)$, which implies that the constraint $\mathcal{K}_2
+\mathcal{K}_3=-T_1^{*3}\e$ is not possible. We only treat case (1) in detail, because cases (2)--(3) follow from 
similar computations. 

For case (1) we need to consider three sub-cases:
\begin{itemize}
\item[(1a)] $\lambda(I)$ is constant and $g_{11}\asymp -\e^{1/3}$.
\item[(1b)] $\lambda(I)^3 = \omega(\e)$, $g_{11}^3=\omega(\e)$ and $\lambda(I)^3g_{11}^3\asymp-\e$.
\item[(1c)] $\lambda(I)^3\asymp\e$ and $g_{11}^3$ is constant and negative.
\end{itemize}
The case $\lambda(I)^3 = o(\e)$ can be excluded, since this would imply that $g_{11}^3\asymp
-\frac{\e}{o(\e)}$, which tends to $-\infty$, a property that is not allowed because $g_{11}\in(-T_1^*,1-T_1^*)$. For 
each of the three sub-cases we study the asymptotic behavior of $\mathcal{K}_2$ as $\e\downarrow0$. 

\paragraph{Case (1a).} 

Since $\lambda(I)$ is constant, it suffices to analyse the square appearing in $\mathcal{K}_2$, i.e., 
\begin{equation}
\left(\frac{\lambda(I)}{1-\lambda(I)}g_{11}-\frac{1-\lambda(I)}{\lambda(I)}g_{22}\right)^2.
\end{equation}
After straightforward computations we see that if $\mathcal{K}_2>0$, then  
\begin{equation}
\label{eq: big big big 2}
\left(\frac{\lambda(I)}{1-\lambda(I)}g_{11}-\frac{1-\lambda(I)}{\lambda(I)}g_{22}\right)^2\asymp \e^{2/3},
\end{equation}
which yields that $\mathcal{K}_2+\mathcal{K}_3 \asymp \e^{2/3}$ instead of $-\e$. Hence this case is not 
possible. 

\paragraph{Case (1b).} 

We have $\lambda(I)=\omega(\e^{1/3})$, $g_{11}=\omega(\e^{1/3})$ and $\lambda(I)^3g_{11}^3
\asymp-\e$, and again obtain that $\lambda(I)^2g_{11}^2\asymp\e^{2/3}$, which yields a result similar to 
the one in \eqref{eq: big big big 2}. 

\paragraph{Case (1c).} 

After straightforward computations we again obtain $\lambda(I)^2g_{11}^2\asymp \e^{2/3}$. 

\medskip
Performing similar computations for cases (2)-(4), we can also exclude those. This verifies the claim 
made in the beginning, namely, that in order to find the optimal graphon corresponding to the constraints 
$\mathcal{K}_1=0$ and $\mathcal{K}_2+\mathcal{K}_3=-T_1^{*3}\e$ it suffices to consider the constraints 
$\mathcal{K}_1=0$, $\mathcal{K}_2=0$ and $\mathcal{K}_3=-T_1^{*3}\e$.


\end{document}